\documentclass[a4paper]{article}
\usepackage[a4paper]{geometry}

\usepackage[cal=esstix,bb=esstix]{mathalpha}
\usepackage[backend=bibtex]{biblatex}
\usepackage[toc,page]{appendix}
\usepackage{hyperref, 
    subfiles,
    amsmath, 
    amsthm, 
    mathtools,
    csquotes,
    diffcoeff,
    microtype,
    enumitem,
    float,
    zref-clever
}

\zcsetup{cap=true}

\theoremstyle{plain}

\theoremstyle{plain}

\AddToHook{env/corollary/begin}{%
\zcsetup{countertype={theorem=corollary}}}

\theoremstyle{plain}

\AddToHook{env/proposition/begin}{%
\zcsetup{countertype={theorem=proposition}}}

\theoremstyle{plain}

\zcRefTypeSetup{lemma}{cap=true}
\AddToHook{env/lemma/begin}{%
\zcsetup{countertype={theorem=lemma}}}

\theoremstyle{remark}
\newtheorem*{remark}{Remark}

\setlist[description]{leftmargin=\parindent,labelindent=\parindent}

\numberwithin{equation}{section}
\zcRefTypeSetup{equation}{
  Name-sg-ab =,
  name-sg-ab =,
  Name-pl-ab =,
  name-pl-ab =,
  abbrev,
}

\usepackage{objects}
\usepackage{main-statements}
\usepackage{supporting-statements}

\title{The Steklov Spectrum of Spherical Cylinders}
\author{Spencer Bullent}
\date{}
\begin{document}

\maketitle

\begin{abstract}
    The Steklov problem on a compact Lipschitz domain is to find harmonic functions on the interior whose outward normal derivative on the boundary is some multiple (eigenvalue) of its trace on the boundary. These eigenvalues form the Steklov spectrum of the domain. This article considers the Steklov spectrum of spherical cylinders (Euclidean ball times interval). It is shown that the spectral counting function admits a two term asymptotic expansion. The coefficient of the second term consists of a contribution from the curvature of the boundary and a contribution from the edges.
\end{abstract}

\section{Introduction and Main Result}\label{section:main-result}
Given a compact domain \(M \subset \mathbb{R}^n\), with Lipschitz boundary \(\partial M\), then the \emph{Steklov problem} is to find the \emph{Steklov spectrum} of eigenvalues, \(\sigma\), with eigenfunctions, \(u_\sigma \in L^2(M)\), satisfying
\begin{equation}
    \label{eqn:SteklovProblem}
    \begin{dcases}
        \Delta u_\sigma = 0 & \text{ in } M\, , \\
        \diffp{u_\sigma}{\nu}  = \sigma u_\sigma & \text{ on } \partial M \, .
    \end{dcases}
\end{equation}
Here, \(\diffp*{}{\nu}\) denotes the outward normal derivative on \(\partial M\), and \(\Delta\) the Laplacian. The nomenclature is justified as the Steklov problem is equivalent to finding the spectrum of the Dirichlet-to-Neumann map on \(M\). From \cite{DtNSemiGroup}, the Steklov spectrum of \(M\) is discrete and accumulates at infinity,
\begin{equation*}
  0 = \sigma_0 \leq \sigma_1 \leq \sigma_2 \leq \cdots \nearrow \infty \,.
\end{equation*}
For pure point spectra that are bounded below the \emph{spectral counting function},
\begin{equation*}
    \NDef \,,
\end{equation*}
is key object that is well studied. It is often sought that the asymptotic expansion of \(N(\sigma)\), as \(\sigma \to \infty\), should follow a \emph{Weyl law}, in that the coefficients depend on topological or geometric properties of the domain. For a general overview and bibliography concerning the Weyl law see \cite{Ivrii2016} and \cite{mårdby2024112yearslisteningriemannian}. Moreover, the Steklov problem has been extensively studied, for an overview and bibliography see \cite{SpectralGeometrySteklov} and \cite{Colbois2024}.

Fix \(n \geq 3,\, \R > 0,\, \oL > 0\). The main result of this article concerns the Weyl law of the Steklov spectrum of \emph{spherical cylinders},
\begin{equation*}
    \Omega \coloneq \Set*{x \in \mathbb{R}^{\n - 1} \given \norm*{x} < \R} \times (-\oL, \oL) \,.
\end{equation*}
\TheoremMain
\begin{remark}
    \item \begin{itemize}
        \item The Steklov spectrum of a spherical cylinder uniquely determines its radius,
        length and dimension.
        \item As, by direct estimation, 
            \begin{equation*}
                \int_{0}^{1} (1 + x)^{-\half}(1 - x)^{\half[\n] - 1} \dl{x} = \int_{0}^{1} \sqrt{\frac{x}{1 + x}} x^{{\half} -1}(1 - x)^{\half[\n] - 1} \dl{x} < \half[\sqrt{2}] B(\thalf, \thalf[\n])\,,
            \end{equation*}
        then, for \(n \geq 3\), \(C_{n, 3} > {\thalf}\frac{\measure*{\mathbb{B}_{\n - 2}}}{(2 \pi)^{\n - 2}}\spar*{2^{\thalf[\n]} - 1 - \sqrt{2}} > 0\). Hence, there is a non-trivial spectral contribution coming from the edges.
        \item For spherical cylinders
        \begin{align*}
            & \measure*{\partial\Omega} = 2 \measure*{\mathbb{B}_{\n - 1}} R^{\n - 1} + 2 L (\n - 1) \measure*{\mathbb{B}_{\n - 1}} R^{\n - 2} \, , \\
            & \measure*{\partial^2 \Omega} = 2 (\n - 1) \measure*{\mathbb{B}_{\n - 1}} R^{\n - 2} \, ,\\
            & \int_{\partial\Omega} \kappa = \frac{\n - 2}{(\n - 1) R} 2 L (\n - 1) \measure*{\mathbb{B}_{\n - 1}} R^{\n - 2} \, .
        \end{align*}
        \item For \(n = 3\) the explicit Weyl law is
        \begin{align*}
            N(\sigma) 
            & = \frac{\pi}{4 \pi^2} (2 \pi R^2 + 4 \pi R L) \sigma^2 + \frac{\pi}{4 \pi^2} (2 \pi L) \sigma + \frac{2}{2 \pi}(4 \pi R) (\sqrt{2} - {\thalf} - (\sqrt{2} - 1)) \sigma + \bigO(\sigma^{\frac{3}{4}}) \\
            & = \spar*{\half[R^2] + R L} \sigma^2 + \spar*{\half[L] + 2 R} \sigma + \bigO(\sigma^{\frac{3}{4}}) \,.
        \end{align*}
        \item The error \(\bigO(\sigma^{\n - 2 - \quarter})\) is not expected to be sharp and could 
        likely be improved to at least \(\bigO(\sigma^{\n - 2 - \third})\) with more detailed analysis.
        \item For asymptotic behaviour in a single variable the notation of Bachmann and Landau is used, in particular \(f(\sigma) = \bigO(\phi(\sigma)) \iff \limsup_{\sigma \to \infty} \Mod*{\frac{f(\sigma)}{\phi(\sigma)}} < \infty\).
    \end{itemize}
\end{remark}

There are a variety of results that \zcref{thm:main} builds upon. The weakest results are for \enquote{singular} boundaries. The \emph{rough} Weyl law,
\begin{equation}
    \label{eqn:lipschitz-counting-function}
    \N(\sigma) = C_{n, 1} \measure*{\partial M} \sigma^{n - 1} + \littleO(\sigma^{n - 1}) \, ,
\end{equation}
has been shown for Lipschitz boundaries in arbitrary dimension (see \cite[Theorem 1.2]{LipschitzAsymptotics}), and for surfaces whose conformal factor, to the circle domain in its conformal class, on the boundary, is in the Orlicz class \(L\log L\) (see \cite[Theorem 1.9]{ConformalAsymptotics}). For more regular boundaries the \emph{sharp} Weyl law
\begin{equation*}
    \N(\sigma) = C_{n, 1} \measure*{\partial M} \sigma^{n - 1} + \bigO(\sigma^{n - 2}) \,,
\end{equation*}
has been shown for surfaces with \(C^{1, 1}\) boundary (see \cite[Theorem 1.12]{CAsymptotics}), in arbitrary dimension for \(C^{2, \alpha}\), \(\alpha > 0\) boundary (see \cite[Theorem 1.11]{CAsymptotics}), and for curvilinear polygons (see \cite[Proposition 2.30]{SloshingSteklovCorners}).

The strongest results are for smooth boundaries. If \(M\) satisfies the dynamical condition, that the set of unit speed periodic geodesics has measure \(0\) in \(S^* M\), then (see \cite{SherHeatInvariants})
\begin{equation}\label{eqn:smooth-counting-function}
    N(\sigma) = C_{n, 1} \measure*{\partial M} \sigma^{n - 1} + C_{n, 2} \spar*{\int_{\partial M} \kappa} \sigma^{n - 2} + \littleO(\sigma^{n - 2}) \, .
\end{equation}
The analysis on smooth boundaries relies on the Dirichlet-to-Neumann map being a classical pseudo-differential operator of order \(1\), whose symbol can be approximated to arbitrary order. Hence, if further terms are known to exist, then the coefficients can, in principle, be calculated from the symbol, and depend on increasingly complex geometric quantities.

A natural class of Lipschitz boundaries to analyse is that of piecewise smooth boundaries, which \(\Omega\) falls into. The motivating example in this case is for cuboids, a product of intervals, in dimension \(n \geq 3\) (see \cite[Theoreom 1.1]{SteklovCuboids})
\begin{equation}
    \label{eqn:cuboid-counting-function}
    \N(\sigma) = C_{n, 1} \measure*{\partial M} \sigma^{n - 1} + C_{n, 3} \measure*{\partial^2 M} \sigma^{n - 2} + \bigO(\sigma^{\eta(n)}) \, .
\end{equation}
\begin{remark}
    The definition of \(C_{n, 3}\) in \cite[Theoreom 1.1]{SteklovCuboids} seems different to that defined above. These two definitions are shown to be equivalent in \zcref{section:edge-equivalence}. 
\end{remark}
\noindent Other related results on piecewise smooth domains include \cite{Ivrii2019} and the sloshing problem (mixed Neumann Robin boundary conditions) in \cite{Levitin2022}.

This article provides another example of a domain with piecewise smooth boundary, with the novelty of the presence of curvature. Comparing to known results, the first term is as expected from \zcref{eqn:lipschitz-counting-function}, whilst the second term has a contribution from the curvature of the smooth components as in \zcref{eqn:smooth-counting-function}, and a contribution from the edges, as if there was no curvature as in \zcref{eqn:cuboid-counting-function}. It is, hence, natural to conjecture that the second term is universal for any piecewise smooth domain whose edges meet at \(\thalf[\pi]\), and that the interaction of curvature and the edges would be visible in lower order terms for smoother spectral aggregates, e.g. the heat trace or Riesz means.

\section{Outline of the Proof of \texorpdfstring{\zcref{thm:main}}{Main Theorem}}\label{section:outline}

The three main steps of the proof are:
\begin{itemize}
    \item Use separation of variables to characterise the eigenvalues as roots of transcendental equations (\zcref{section:Eigenfunctions-Eigenvalues}). Partition the eigenvalues into \enquote{transversely localised} and \enquote{radially localised} sets.
    \item Express the count of transversely localised (\zcref{section:lattice-point-bessel}) and radially localised (\zcref{section:lattice-point-mod-bessel}) eigenvalues as a weighted lattice point count in some domain and approximate this domain.
    \item Approximate the weighted lattice point counts using Euler-Maclaurin and Van der Corput's method (\zcref{section:lattice-point-count-approx}).
\end{itemize}

\begin{remark}
By radially localised it is meant that the eigenfunctions exponentially decay as \(r \to 0\), and exhibit wave like structure in the remaining dimensions, i.e. waves are localised on the \enquote{sleeve} of the spherical cylinder. Similarly, transversely localised means that the eigenfunctions exponentially decay as \(z \to 0\) and exhibit wave like structure in the remaining dimensions, i.e. waves are localised to the \enquote{ends} of the spherical cylinder.
\end{remark}

In more detail, note that \(\partial \Omega\) is piecewise smooth, with smooth components denoted
\begin{equation*}
    \GammaDef,\qquad \GammaPmDef \,.
\end{equation*}
In cylindrical co-ordinates,
\begin{equation*}
    (r, x, z) \in [0, \R] \times S^{\n - 2} \times [\oL, \oL] \mapsto (r x, z) \in \overline{\Omega} \, ,
\end{equation*}
the Steklov problem on \(\Omega\) is
\begin{align}\label{eqn:radial_problem}
    \begin{dcases}
        r^{-(n - 2)} \difcp{r^{n - 2} \difcp{u}{r}}{r} + r^{-2} \Delta_{S^{n - 2}} u + \difcp[2]{u}{z} = 0 & \text{ in } \Omega \, ,\\
        \eval{\difcp{u}{r}}_{r = \R} = \sigma u(R, x, z) & \text{ on } \Gamma \, ,\\
        \eval{\pm \difcp{u}{z}}_{z = \pm \oL} = \sigma u(r, x, \pm \oL) & \text{ on } \Gamma^\pm \, .
    \end{dcases}
\end{align}
Here \(S^{n}\) is the unit sphere in \(\mathbb{R}^{n + 1}\) and \(\Delta_{S^{n}}\) is the Laplace-Beltrami operator associated to the round metric.

The spectral problem \zcref{eqn:radial_problem} admits solutions of separated variables, and these solutions give a complete description of eigenvalues. In particular, denote
\begin{itemize}
    \item \(\Jv\) as the principal branch of the Bessel function of the first kind, defined in \zcref{eqn:bessel-j-def},
    \item \(\Iv\) as the principal branch of the modified Bessel function of the first kind, defined in \zcref{eqn:bessel-i-def},
    \item \(\omega_{n, k}\) as the dimension of the eigenspace of the spherical harmonics on \(S^n\) of order \(k\), defined in \zcref{eqn:spherical-multiplicity}.
\end{itemize}
\PropProblemSolutions[Proved in \zcref{section:Eigenfunctions-Eigenvalues}]

Using \zcref{prop:problem-solns} the spectral counting function may be written as
\begin{equation}
    \label{eqn:structured-counting-func}
    \StructuredCountingFunc \, .
\end{equation}
Consider the asymptotic behaviour of \(\N(\sigma)\) as \(\sigma \to \infty\). There are \(5\) sums over \(k\) to analyse, the simplest one to consider first is the sum over the exceptional eigenvalue, which takes values of either \(0\) or \(\omega_{n - 2, k}\),
\begin{equation*}
    \sum_{k = 0}^{\infty} \omega_{n - 2, k} \mathbb{1}_{\Set*{\oL \R = k, \oL < \sigma}} = \bigO(1) \, .
\end{equation*}

Next, the sums containing \(\Sigma_{1, k}\) and \(\Sigma_{2, k}\), corresponding to transversely localised eigenvalues, are similar to each other, and are analysed together. For intuition consider the transcendental equations
\begin{align*}
    & x \frac{\Jv_1'(x)}{\Jv_1(x)} = x \coth(\tquarter[x]) \, , \\
    & x \frac{\Jv_1'(x)}{\Jv_1(x)} = x \tanh(\tquarter[x]) \, .
\end{align*}
Graphically
\begin{figure}[H]
    \center{
        \includegraphics{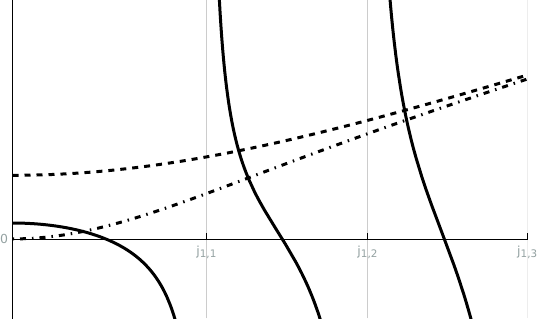}
        \caption{
            \label{fig:example-graph}
            \(x \frac{\Jv_1'(x)}{\Jv_1(x)}\) - solid, \(x \coth(\tquarter[x])\) - dashed, \(x \tanh(\tquarter[x])\) - dot dashed 
        }
    }
\end{figure}
\noindent Where \(j_{\nu, k}\) denotes the \(k\)-th positive root of \(\Jv\). Informally \zcref{fig:example-graph} indicates that \(x \frac{\Jv_1'(x)}{\Jv_1(x)}\) is monotonically decreasing on intervals of continuity, has multi-valued inverse, and hence will have an intersection with the monotonically increasing \(x \coth(\tquarter[x])\) and \(x \tanh(\tquarter[x])\) in the ranges \((j_{\nu, j}, j_{\nu, j + 1})\) and a potential intersection in \((0, j_{\nu, 1})\). This can be made rigorous and allows both sums to be approximated by counting the positive lattice points in a function's subgraph. In particular, denote
\begin{itemize}
    \item \(\lfloor x \rfloor\) as the integer part of \(x\),
    \item \(\Ai\) as the standard Airy function, defined in \zcref{eqn:ai-def},
    \item \(a_k\) the \(k\)-th root of \(\Ai\),
    \item \(\Yv\) as the principal branch of the Bessel function of the second kind, defined in \zcref{eqn:bessel-y-def}
    \item \(m(\nu,x), \theta(\nu, x)\) as the \emph{Bessel-Riccati modulus} and \emph{phase} respectively, defined in \zcref{eqn:def_bessel_modulus} and \zcref{eqn:def_bessel_phase}, informally \begin{equation*}
        \Jv(x) + \iu \Yv(x) = \sqrt{\frac{2}{\pi x}} m(\nu, x) e^{\iu \theta(\nu, x)} \, .
    \end{equation*}
\end{itemize}
Moreover, define
\begin{align}
    & \BetaDef \, , \nonumber \\
    & \label{eqn:psidef} \begin{aligned}
        & \PsiDomain \\
        & \PsiDef \, ,
    \end{aligned} \\
    & \xtDef \, , \nonumber \\
    & \xcDef \, . \nonumber
\end{align}
\PropBesselSumAsymptotics[Proved in \zcref{section:lattice-point-bessel}]

Furthermore, in order to give accurate asymptotics for the sums in \zcref{prop:bessel-sum-asymptotics} one requires a controlled approximation of the summand. In particular, define
\begin{equation*}
    \begin{aligned}
        & \EtaOneDomain \\
        & \EtaOneDef \, .
    \end{aligned}
\end{equation*}
\PropBesselApprox[Proved in \zcref{section:lattice-point-bessel}]

Next, the sums containing \(\Sigma_{3, k}\) and \(\Sigma_{4, k}\), corresponding to radially localised eigenvalues, can be analysed together and be expressed as a lattice point count. For \(\nu \geq 0\) define
\begin{equation*}
    \xvDef \, .
\end{equation*}
\ModBesselSumAsymptotics[Proved in \zcref{section:lattice-point-mod-bessel}]

Again, the sums in \zcref{prop:mod-bessel-sum-asymptotics} require a controlled approximation of the summand. It is sufficient to approximate
\begin{equation*}
    \frac{\oL}{\R} \xv_{\nu}(\sigma) + \fcall{\arctan}{\frac{\sigma - \beta}{\xv_{\nu}(\sigma)}} \, ,
\end{equation*}
with \(\nu\) replaced by \(k + \beta\), and \(\sigma\) replaces with \(\sigma \R + \beta\).
Define
\begin{equation*}
    \begin{aligned}
        & \EtaTwoDomain \\
        & \EtaTwoDef \, .
    \end{aligned}
\end{equation*}
\PropModBesselApprox[Proved in \zcref{section:lattice-point-mod-bessel}]

Finally, the lattice point count requires approximation. Informally, given a \enquote{nice} enough function \(g\), which in practice will be \(\omega_{n - 2, k}\spar*{\frac{\eta_1}{\pi} + c}\) or \(\omega_{n - 2, k}\spar*{\frac{\eta_2}{\pi} + c}\), then 
\begin{equation*}
    \sum_{j = 0}^{n} \Floor{g(j)} = \sum_{j = 0}^{n} g(j) + \varepsilon_1 = \int_{0}^{n} g(x) \dl{x} + \varepsilon_1 + \varepsilon_2 \, .
\end{equation*}

The first error, \(\varepsilon_1\), can be estimated via Van der Corput's method. Denote \(\rho_1\) as the first Euler-Maclaurin error function, otherwise known as the \emph{rounding error} or \emph{sawtooth} function,
\begin{equation*}
    \EMErrorOneDef
\end{equation*}
\PropVanDerCorputErrorEstimate

The second error, \(\varepsilon_2\), can be estimated with the Euler-Maclaurin formula (see e.g. \cite[\S 8.1]{olver1997asymptotics}). Given a function \(f \in C^1([a, b])\), then
\begin{equation}
    \label{eqn:euler-maclaurin}
    \EMOne \, ,
\end{equation}
moreover, denoting \(\mathcal{V}_{a}^{b}(f)\) as the variation of \(f\) over \([a, b]\), then
\begin{equation*}
    \EMOneRemainderBound \, .
\end{equation*}

In particular, the following estimates hold for the approximation functions \(\eta_1, \eta_2\)
\PropBesselStrongConvexity[Proved in \zcref{section:lattice-point-count-approx}]
\PropModifiedBesselStrongConcavity[Proved in \zcref{section:lattice-point-count-approx}]

\zcref{prop:BesselStrongConvexity} and \zcref{prop:ModifiedBesselStrongConcavity} can be applied to the weighted sum via the following theorem.
\TheoremHarmonicWeightSlicing

After applying \zcref{thm:HarmonicWeightSlicing} it remains to approximate various integrals. For the level of approximation sought there is a single non-elementary term within \(C_{n, 3}\) that is denoted
\begin{equation}\label{eqn:G-edge-constant}
    G'_{n, 1} \coloneq \frac{\int_{0}^{1} (1 + x)^{-\half}(1 - x)^{\half[n] - 1} \dl{x}}{B(\thalf, \thalf[n])} \, ,
\end{equation}
in keeping with the definition of \(G_{d - 1, 1}\) in \cite[Theorem 1.1]{SteklovCuboids}.

As an example of putting all the preceding statements together, consider the sum over \(\Sigma_{3, k}\) for \(\oL = \R = 1\). The following steps can be taken, in an informal sense, in that no bookkeeping of the error is performed.
\begin{enumerate}
    \item Express the sum over \(\Sigma_{3, k}\) in terms of a sum over tractable functions using \zcref{prop:mod-bessel-sum-asymptotics}
    \begin{flalign*}
        & \sum_{k = 0}^{\infty} \omega_{n - 2, k} \#\Set{\sigma_{3, k, \ell} \in \Sigma_{3, k} \given \sigma_{3, k, \ell} < \sigma} &\\
        & \qquad \sim \sum_{k = 0}^{\sigma} \omega_{n - 2, k} \Floor*{\frac{\xv_{k + \beta}(\sigma + \beta) + \fcall{\arctan}{\frac{\sigma}{\xv_{k + \beta}(\sigma + \beta)}}}{\pi}} &
    \end{flalign*}
    \item Approximate the summands with elementary functions using \zcref{prop:mod-bessel-lattice-point-func-approx}
    \begin{flalign*}
        & \qquad \sim \sum_{k = 0}^{\sigma} \omega_{n - 2, k} \Floor*{\frac{\eta_2(k + \beta, \sigma + \beta) + \thalf[\pi]}{\pi}} &
    \end{flalign*}
    \item Re-weight the sum using \zcref{thm:HarmonicWeightSlicing} into a nested sum
    \begin{flalign*}
        & \qquad \sim \sum_{k = 0}^{\sigma} \binom{k + n - 4}{n - 4} \sum_{j = 0}^{\sigma - k} \omega_{1, j} \Floor*{\frac{\eta_2(j + \beta, \sigma + \beta) + \thalf[\pi]}{\pi}} &
    \end{flalign*}
    \item Approximate the inner sums by an integral using \zcref{prop:ModifiedBesselStrongConcavity}
    \begin{flalign*}
        & \qquad \sim \sum_{k = 0}^{\sigma} \binom{k + n - 4}{n - 4} 2 \int_k^\sigma \frac{\eta_2(\nu + \beta, \sigma + \beta)}{\pi} \dl{\nu} &
    \end{flalign*}
    \item Interchange summation and integration and then approximate the resulting integral
    \begin{flalign*}
        & \qquad \sim \frac{\Mod*{\mathbb{B}_{n - 1}}}{(2 \pi)^{n - 1}} \Mod*{S^{n - 1}} \sigma^{n - 1} + \spar*{\frac{\Mod*{\mathbb{B}_{n - 1}}}{(2 \pi)^{n - 1}} \half[(n - 2)(n - 1)] \Mod*{S^{n - 1}} - \frac{\Mod*{\mathbb{B}_{n - 2}}}{(2 \pi)^{n - 2}} \frac{G'_{n - 1, 1}}{4}} \sigma^{n - 2} &
    \end{flalign*}
\end{enumerate}

\zcref{thm:main} then follows from similar analysis on the remaining sums, with the appropriate propositions and calculations used.

\subsection{Further Notation and Definitions}
Further notation not stated above, and some basic definitions, are collected here.
\begin{itemize}
    \item \(\arctan(x) \coloneq \int_{0}^{x} \frac{\dl{y}}{1 + y^2}\) denotes the principal branch of \(\arctan\) restricted to \(\mathbb{R}\),
    \item \(W(f, g)\) denotes the Wronskian of two differentiable functions,
    \item \(\BetaFDef\), \(\IncompleteBetaDef\) denote the Beta function and incomplete Beta function respectively,
    \item \(\PochammerDef\) denotes the Pochhammer symbol,
    \item \(\Bi\) denotes the \(2\)-nd standard Airy function, defined in \zcref{eqn:bi-def},
    \item \(\AiryModulusDef\),
    \item \(\Hv^1, \Hv^2\) denote the principal branch of the Bessel functions of the third kind, defined in \zcref{eqn:bessel-3-def}
    \item \(j'_{\nu, k}, y_{\nu, k}\) denote the \(k\)-th positive root of \(\Jv'(x), \Yv(x)\) respectively,
    \item \(\tilde{\theta}(\nu, x), \varepsilon_\theta(\nu, x)\) denote an upper approximation to the Bessel phase \(\theta\), and an upper bound on the error of \(\tilde{\theta}\) respectively, defined in \zcref{eqn:def_bessel_modulus}, \zcref{eqn:def_bessel_phase},
    \item \(\xi(\nu, x)\) denotes the variable used in Olver's Bessel function asymptotics, defined in \zcref{eqn:xiDef},
    \item \(K_\nu\) denotes the principal branch of the modified Bessel function of the second kind,
    \item \(p(z), \zeta(z), U_1(p), V_1(p)\) denote variables and functions used in Olver's modified Bessel asymptotics, defined in \zcref{eqn:zetaDef},
    \item \(\phi(\nu, x)\) denotes the phase between \(x \Jv'(x)\) and \(\Jv(x)\), defined in \zcref{prop:BesselConditionBehaviour},
    \item \(I(z; p, q)\) denotes part of \(C_{n, 3}\) defined in \zcref{eqn:edge-integral}.
\end{itemize}

\section{Steklov Eigenfunctions and Transcendental Equations Governing the Associated Eigenvalues}\label{section:Eigenfunctions-Eigenvalues}
Recall that the first step is \zcref{prop:problem-solns}, which shows the structure of the Steklov eigenvalues, namely as collections of sets depending on transcendental equations. The following definitions are required:
\begin{itemize}
    \item Define \(\Jv(x) \in \BesselFunctionSpace\) as the solution to the differential equation
    \begin{equation}
        \label{eqn:bessel-j-def}
        \begin{dcases}
            \BesselDifferentialEquation[\Jv][\nu][x] & \mid x > 0 \, ,\\
            \JvBC[\nu][x] \, .
        \end{dcases}
    \end{equation}

    \item Define \(\Iv(x) \in \ModifiedBesselFunctionSpace\) as the solution to the differential equation
    \begin{equation}
        \label{eqn:bessel-i-def}
        \begin{dcases}
            \ModifiedBesselDifferentialEquation[\Iv][\nu][x] & \mid x > 0 \, ,\\
            \IvBC[\nu][x] \, .
        \end{dcases}
    \end{equation}

    \item The spherical harmonics on \(S^{n}\) of order \(k \in \mathbb{N}_{0}\) are defined as an orthonormal basis to the eigenspace
    \begin{equation*}
        \begin{dcases}
            \SphericalHarmonicDifferentialEquation \, ,\\
            \norm*{y}_{L^2} = 1 \, ,
        \end{dcases}
    \end{equation*}
    where each eigenvalue \(k(k + n - 1)\) has multiplicity
    \begin{equation}
        \label{eqn:spherical-multiplicity}
        \SphericalHarmonicMultiplicityDef
    \end{equation}
\end{itemize}
Moreover, the following Lemma will be used.
\LemmaSeperableEigenbasis

\begin{proof}[Proof of \zcref{prop:problem-solns}]
    The Steklov problem in cylindrical co-ordinates, \zcref{eqn:radial_problem}, admits separable solutions, hence, consider the ansatz,
    \begin{equation*}
        u(r, x, z) = u_r(r) u_x(x) u_z(z) \, .
    \end{equation*}
    First, in the spherical variable \(x\), consider finding \(u_x \in L^2(S^{n - 2})\) such that 
    \begin{equation*}
            \Delta_{S^{n - 2}} u_x = - \lambda u_x \,.
    \end{equation*}
    This spectral problem has solutions of linear combinations of spherical harmonics on \(S^{n - 2}\) of order \(k \in \mathbb{N}_{0}\) and associated eigenvalue \(\lambda = k (k + n - 3)\). Next, in the longitudinal variable \(z\) consider finding \(u_z \in L^2([-L, L])\) such that 
    \begin{equation*}
            \diff[2]{u_z}{z} = \pm \alpha^2 u_z \,.
    \end{equation*}
    This spectral problem has solutions of linear combinations of
    \begin{equation*}
        \begin{cases}
            1,\, z & \mid \sign(\pm \alpha^2) = 0 \,,\\
            \cosh(\alpha z),\, \sinh(\alpha z) & \mid \sign(\pm \alpha^2) = 1 \,,\\
            \cos(\alpha z),\, \sin(\alpha z) & \mid \sign(\pm \alpha^2) = -1 \,.
        \end{cases}
    \end{equation*}
    Lastly, in the radial variable \(r\), consider finding \(u_r \in L^2([0, R])\) such that
    \begin{equation*}
            r^{-(n - 2)}\diff*{\spar*{r^{n - 2}\diff{u_r}{r}}}{r} + \spar*{\pm \alpha^2 - \frac{k(k + n - 3)}{r^2}}u_r = 0 \,.
    \end{equation*}
    This differential equation can be rewritten as
    \begin{equation*}
        \diff*[2]{(r^{\frac{n}{2} - 1}u_r(r))}{r} + \spar*{\pm \alpha^2 - \frac{(k + \frac{n}{2} - \frac{3}{2})^2 - \frac{1}{4}}{r^2}} r^{\frac{n}{2} - 1} u_r(r) = 0 \, .
    \end{equation*}
    Hence, the spectral problem has known solutions of
    \begin{equation*}
        r^{\frac{n}{2} - 1} u_r(r) = \begin{cases}
            r^{k + \frac{n}{2} - 1} & \mid \sign(\pm \alpha^2) = 0 \,,\\
            r^{\frac{1}{2}} \Jv_{k + \frac{n}{2} - \frac{3}{2}}(\alpha r) & \mid \sign(\pm \alpha^2) = 1 \,,\\
            r^{\frac{1}{2}} \Iv_{k + \frac{n}{2} - \frac{3}{2}}(\alpha r) & \mid \sign(\pm \alpha^2) = -1 \,.
        \end{cases}
    \end{equation*}
    Next, considering the boundary conditions
    \begin{align*}
        \begin{dcases}
            \eval{\diff{u_r(r)}{r}}_{r = R} = \sigma u_r(R) \,,\\
            \eval{\diff{u_z(z)}{z}}_{z = L} = \sigma u_z(L) \,,\\
            \eval{-\diff{u_z(z)}{z}}_{z = - L} = \sigma u_z(- L) \,,
        \end{dcases}
    \end{align*}
    then \(u_z, \diff{u_z}{z}\) is either even and odd respectively, or vice-versa and so must be a multiple of \(z, \cosh, \sinh, \cos,  \sin\), depending on \(\sign(\pm \alpha^2)\). Moreover, equating \(\sigma\) leads to the required transcendental equation for \(\alpha\). Lastly by \zcref{lem:seperable-eigenbasis} every eigenspace of \zcref{eqn:radial_problem} has a basis of separable solutions and so the above covers all eigenvalues and multiplicities, upon choosing a basis of the spherical harmonics for each \(k \in \mathbb{N}_{0}\). 
Hence, the spectral counting function can be expressed in a structured way
\begin{equation*}
    \StructuredCountingFunc \, .    
\end{equation*}
\end{proof}

\section{Collected properties of Airy and Bessel functions}
Analysing \(\N(\sigma)\) requires many properties of both the Airy and Bessel functions. These properties are collected below and can be referenced as required.

\subsection{Properties of Airy functions}

The standard solutions \(\Ai, \Bi, \Ai(e^{\pm \iu \third[2 \pi]} z)\) of the Airy differential 
equation,
\begin{equation*}
    \AiryDifferentialEquation \, ,
\end{equation*} are all entire functions and can be defined by
\begin{align}
    \label{eqn:ai-def}
    & \AiDef \, , \\
    \label{eqn:bi-def}
    & \AiryConnection \, .
\end{align}

The Airy modulus function is defined as \(\AiryModulusDef[z]\) and, for \(x > 0\), has the following envelope \cite[\S 11.1, 11.2]{olver1997asymptotics}
\begin{align}
    \label{eqn:AiryModulusAsymptotic}
    & 1  - \frac{5}{32} x^{-3} < \pi \sqrt{x} M_A^2(-x) < 1 \\
    \label{eqn:DAiryModulusAysmptotic}
    & {\half} - \frac{35}{64}x^{-3} < \pi x^{\half[3]} \spar*{M_A^2}'(-x) < {\half}
\end{align}

\subsection{Properties of Bessel Functions}
The additional standard solutions, \(\Yv, \Hv^1, \Hv^2\), of the Bessel differential equation,
\begin{equation}
    \label{eqn:Bessel}
    \BesselDifferentialEquation\, ,
\end{equation}
can be defined via the asymptotic behaviour as \(x \to \infty\)
\begin{equation}
    \label{eqn:bessel-y-def}
    \YvBC \, ,
\end{equation}
and the connection formulae
\begin{align}
    \label{eqn:bessel-3-def}
    & \BesselConnectionI[\nu][x] \, ,\\
    & \BesselConnectionII[\nu][x] \, .
\end{align}

\zcref{eqn:Bessel} has a singularity as \(x \to 0\), leading to the following asymptotics for \(\nu \geq 0\) fixed, as \(x \to 0\) \cite[\S 3.1 eq. (8) \S 3.54 eq. (1), (2)]{WatsonBessel}
\begin{align}
    \label{eqn:BesselSmallArg}
    & \Jv(x) = x^\nu \spar*{\frac{1}{2^\nu \Gamma(\nu + 1)} + \bigO(x^2)} \, ,\\
    \label{eqn:DBesselSmallArg}
    & \Jv'(x) = x^{\nu - 1} \spar*{\frac{\nu}{2^\nu \Gamma(\nu + 1)} + \bigO(x^2)} \, ,\\
    \label{eqn:Bessel2SmallArg}
    & \Yv(x) = \begin{dcases}
        \frac{2}{\pi}\log(\thalf[x]) + \bigO(1) & \mid \nu = 0 \, ,\\
        (\thalf[x])^{-\nu}\spar*{-\frac{\Gamma(\nu)}{\pi} + \bigO(x^{2 \nu} + x^2)} & \mid \nu > 0 \, .
    \end{dcases}
\end{align}

\zcref{eqn:Bessel} also has a turning point at \(x = \nu\), after which the Bessel functions are oscillatory. From \cite[\S 15.22]{WatsonBessel} the roots, \(j_{\nu, k}\), \(j'_{\nu, k}\) and \(y_{\nu, k}\), of \(\Jv\), \(\Jv'\) and \(\Yv\) interlace for \(\nu > -1\),
\begin{equation}\label{eqn:BesselRootInterlacing}
    \nu < j'_{\nu, 1} < y_{\nu, 1} < j_{\nu, 1} < \cdots < j'_{\nu, k} < y_{\nu, k} < j_{\nu, k} < \cdots \nearrow \infty \, .
\end{equation}
Moreover, \(j_{\nu, k}\) has bounds, for \(\nu > 0, k \geq 1\), \cite[eq. 1.2]{QuBestPossible}
\begin{equation}
    \label{eqn:QuBestPossible}
    \nu - a_k \spar*{\frac{\nu}{2}}^{\third} < j_{\nu, k} < \nu - a_k \spar*{\frac{\nu}{2}}^{\third} + \frac{3}{20}a_k^2 \spar*{\frac{\nu}{2}}^{-\third} \, .
\end{equation}

Next, consider Kummer's modulus-phase formulation for \(\Jv\) and \(\Yv\). From \zcref{eqn:BesselSmallArg} and \zcref{eqn:Bessel2SmallArg}
\begin{align*}
    W(\Jv(x), \Yv(x))  & = \frac{2}{\pi x} \, ,\\
    \fcall{\arctan}{\frac{\Yv}{\Jv}}(0+) & = - \half[\pi] \, ,
\end{align*}
and from \zcref{eqn:BesselRootInterlacing} \(\Jv\) and \(\Yv\) have no common zeroes. Hence,
\begin{align}
    \label{eqn:def_bessel_modulus}
    & \begin{aligned}
        & \BesselModulusDomain \\
        & \BesselModulusDef \,,
    \end{aligned} \\
    \shortintertext{and}
    \label{eqn:def_bessel_phase}
    & \begin{aligned}
        & \BesselPhaseDomain \\
        & \BesselPhaseDef \,,
    \end{aligned}
\end{align}
are well-defined, with
\begin{equation*}
    \sqrt{\frac{2}{\pi x}} m(\nu, x) e^{\iu \theta(\nu, x)} = \Jv(x) + \iu \Yv(x)\, ,
\end{equation*}
\(\theta(\nu, x)\) is strictly increasing in \(x\) and, given \(x > 0\), \(x \neq j_{\nu, k}\),
\begin{equation*}
    \theta(\nu, x) = \#\{j_{\nu, k} < x\} \pi + \fcall{\arctan}{\frac{\Yv(x)}{\Jv(x)}} \, .
\end{equation*}

Furthermore, \zcref{eqn:Bessel} has a singularity as \(\nu \to \infty, x \to \infty\) with complicated behaviour depending on how the singularities are approached (see \cite{Sher:BesselAsymptotics} for details). Starting with the order dominant part of the Debye scale. From \cite[\S 8.4]{WatsonBessel}, for a compact \(Z \subset (0, 1)\) there exists a constant \(C\) such that, for \(\nu > C\) and \(z \in Z\),
\begin{align}
    \label{eqn:BesselDebyeAsymptotic}
    & \Mod*{\Jv(\nu z) \frac{\sqrt{2 \pi \nu} (1 - z^2)^{\quarter}}{e^{-\nu ((1 - z^2)^{\half} - \arcsech(z))}} - \spar*{1 + \frac{1}{24 \nu} \spar*{3(1 - z^2)^{-\half} - 5 (1 - z^2)^{-\half[3]}}}} \leq C \nu^{-2} \,,\\
    \label{eqn:DBesselDebyeAsymptotic}
    & \begin{aligned}[b]
        & \Biggl\lvert \Jv'(\nu z) \frac{\sqrt{2 \pi \nu} z}{e^{-\nu ((1 - z^2)^{\half} - \arcsech(z))}}(1 - z^2)^{-\quarter} \\ 
        & \hspace{4.4cm} - \spar*{1 + \frac{1}{24 \nu}\spar*{-9(1 - z^2)^{-\half} + 7 (1 - z^2)^{-\half[3]}}} \Biggr\rvert
    \end{aligned} \leq C \nu^{-2} \, .
\end{align}
Next, considering the region \(0 < \nu < x\), define
\begin{equation}
    \label{eqn:xiDef}
    \begin{aligned}
        & \XiDomain \\
        & \XiDef \, .
    \end{aligned}
\end{equation}
From \cite[\S 11.10]{olver1997asymptotics} there exists a constant \(C > 0\) such that, for \(x > C\) and \(0 < \nu < x\),
\begin{align}
    \label{eqn:HankelUniformAsymptotic}
    & \Mod*{\Hv^{j}(x) \nu^{\third} \spar*{\frac{x^2 - \nu^2}{4 \xi \nu^2}}^{\quarter} \spar*{\fcall{\Ai}{- \nu^{\third[2]} \xi} \pm \iu \fcall{\Bi}{- \nu^{\third[2]} \xi}}^{-1} - 1} < C \nu^{-1}\min(1, \xi^{-\half[3]}) \, , \\
    \label{eqn:DHankelUniformAsymptotic}
    & \begin{aligned}
        \Mod*{\spar*{\Hv^j}'(x) \frac{x}{2 \nu^{\third}}\spar*{\frac{4 \xi \nu^2}{x^2 - \nu^2}}^{\quarter} \spar*{\Ai'(- \nu^{\third[2]} \xi) \pm \iu \Bi'(- \nu^{\third[2]} \xi)}^{-1} - 1} < C \nu^{-1}\min(1, \xi^{-\half[3]}) \, ,
    \end{aligned}
\end{align}
where \(j = 1\) corresponds to \(+\) and \(j = 2\) corresponds to \(-\).
Moreover, define
\begin{align*}
    &  \BesselPhaseApproxDomain \\
    &  \BesselPhaseApproxDef \, , \\
    &  \BesselPhaseApproxErrorDomain \\
    &  \BesselPhaseApproxErrorDef \, .
\end{align*}
From \cite[\S 1.4]{filonov2024uniform}, \(\theta(\nu, x)\) has the following envelope
\begin{equation}\label{eqn:BesselPhaseBounds}
    \widetilde{\theta}(\nu, x) - \varepsilon_\theta(\nu, x) \leq \theta(\nu, x) \leq \widetilde{\theta}(\nu, x) \, .
\end{equation}

\begin{remark}
    The above formulation is a trivial improvement over \cite[\S 1.4]{filonov2024uniform} by noting that \(\difcp{(\widetilde{\theta} - \theta)}{x} < 0\) \cite[\S 13.74 eq. (1)]{WatsonBessel} and so
    \begin{align*}
        \sup_{x \geq \nu}  (\widetilde{\theta}(\nu, x) - \theta(\nu, x)) = \widetilde{\theta}(\nu,\nu) - \theta(\nu, \nu) < \quarter[\pi] \, .
    \end{align*}
\end{remark}

\subsection{Properties of Modified Bessel Function}
The defining differential equation of \(\Iv\),
\begin{equation*}
    \ModifiedBesselDifferentialEquation \, ,
\end{equation*}
has two singularities, at \(x \to 0+\) and at \(\nu \to \infty, x \to \infty\), though the behaviour at \(\infty\) is simpler than for Bessel functions.

First, from \cite[eq. 6, 16, 17]{amos:modifiedbessel}, for \(x > 0\) and \(\nu\) fixed, \(\nu \geq 0\), 
\begin{align}
    \label{eqn:ModifiedConditionFuncBound}
    & \frac{x \Iv'(x)}{\Iv(x)} \leq \sqrt{x^2 + \nu^2} \, , \\
    \label{eqn:LimitModifiedConditionFunc}
    & \lim_{x \to 0+} \frac{x \Iv'(x)}{\Iv(x)} = \nu \, , \\
    \label{eqn:MonotonicityModifiedConditionFunc}
    & \frac{x \Iv'(x)}{\Iv(x)} \nearrow \infty \text{ as } x \to \infty \, .
\end{align}

Last, consider the region \(0 < \nu < x\), define
\begin{equation}
    \label{eqn:zetaDef}
    \begin{aligned}
        & z \coloneq \frac{x}{\nu} \, , \\
        & \pDef \, ,\\
        & \zetaDef \, ,\\
        & \UIDef \, ,\\
        & \VIDef \, .
    \end{aligned}
\end{equation}
From \cite[\S 10.7]{olver1997asymptotics} there exists a constant \(C\) such that for \(\sqrt{x^2 + \nu^2} > C\)
\begin{align}
    \label{eqn:ModifiedBesselAsymptotics}
    & \Mod*{\Iv(x) \frac{(2 \pi \nu)^{\half} (1 + z^2)^{\quarter}}{\fcall{\Re}{e^{-\nu \zeta(z e^{\iu \pi})}}} - \spar*{1 + \frac{U_1(p)}{\nu}}} < C \spar*{x^2 + \nu^2}^{-1} \, , \\
    \label{eqn:DModifiedBesselAsymptotics}
    & \Mod*{\Iv'(x) \frac{(2 \pi \nu)^{\half} z}{\fcall{\Re}{e^{-\nu \zeta(z e^{\iu \pi})}} (1 + z^2)^{\quarter}} - \spar*{1 + \frac{V_1(p)}{\nu}}} < C \spar*{x^2 + \nu^2}^{-1} \, .
\end{align}

\section{Lattice Point Count Formulation Over Transversely Localised Waves}\label{section:lattice-point-bessel}
The next step is to consider sums containing \(\Sigma_{1, k}\) and \(\Sigma_{2, k}\). These correspond to transversely localised waves, in that the eigenfunctions exponentially decay as \(z \to 0\). From \zcref{prop:problem-solns} the associated transcendental equations are of the form:
\begin{align*}
    & \alpha \frac{\Jv_{k + \beta}'(\alpha \R)}{\Jv_{k + \beta}(\alpha \R)} - \frac{\beta}{\R} = \alpha \tanh(\alpha \oL)\, , \\
    & \alpha \frac{\Jv_{k + \beta}'(\alpha \R)}{\Jv_{k + \beta}(\alpha \R)} - \frac{\beta}{\R} = \alpha \coth(\alpha \oL) \, .
\end{align*}
Both equations have the structure of a piecewise smooth function with multi-valued inverse equaling a monotonically increasing smooth function. The plan of attack is to mimic
\begin{equation*}
    x = \tan(x) \iff \exists\, n \in \mathbb{Z} \text{ s.t. } \arctan(x) - x = n \pi \, .
\end{equation*}
To start, consider where the possible roots can be. Note that the first transcendental equation can be rewritten as
\begin{equation*}
    \alpha \R \frac{\Jv'(\alpha R)}{\Jv(\alpha R)} = \beta + \alpha \R \fcall{\tanh}{\alpha \R \frac{\oL}{\R}} \, ,
\end{equation*}
and so, it is sufficient to consider the roots of 
\begin{equation}
    \label{eqn:BesselRoots}
    \frac{x \Jv'(x)}{\Jv(x)} = f(x) \, ,
\end{equation}
for some continuous, monotonic and increasing function \(f\). In application \(f\) will be either \(\beta + x \tanh(x \oL / \R)\) or \(\beta + x \coth(x \oL / \R)\). The multi-valued nature of \(\Yv / \Jv\) is captured in the phase \(\theta\), and so, it is reasonable to define the same for \( x \Jv' / \Jv\).
\PropBesselConditionBehaviour
\begin{proof}
    As \(x \to 0\), using \zcref{eqn:BesselSmallArg} and \zcref{eqn:DBesselSmallArg},
    \begin{equation*}
    \fcall{\arctan}{\frac{x \Jv'(x)}{\Jv(x)}} = \fcall{\arctan}{\nu \frac{c + \bigO(x^2)}{c + \bigO(x^2)}} = \arctan(\nu) + \bigO(x^2) \, .
    \end{equation*}
    Similarly,
    \begin{align*}
    \diff{}{x}\fcall{\arctan}{\frac{x \Jv'(x)}{\Jv(x)}} 
        & = \frac{W(\Jv(x), x \Jv'(x), )}{\Jv^2(x) + (x \Jv'(x))^2} \\
        &= \frac{(\nu^2 - x^2) \Jv^2(x) - (x \Jv'(x))^2}{x (\Jv^2(x) + (x \Jv'(x))^2)} \\
        &= -\frac{x^2 (c^2 + \bigO(x^2))}{x(c^2(1 + \nu^2) + \bigO(x^2))} \\
        &= \bigO(x)
    \end{align*}
    and
    \begin{align*}
    \diff{}{x} \spar*{x W(\Jv(x), x \Jv'(x))} 
        & = \diff{}{x} \spar*{(\nu^2 - x^2) \Jv^2(x) - (x \Jv'(x))^2} \\
        & = - 2 \spar*{x \Jv^2(x) + \Jv'(x) \spar*{x^2 \Jv''(x) + x \Jv'(x) + (x^2 - \nu^2) \Jv(x)}} \\
        & = -2 x \Jv(x)^2 \leq 0 \, .
    \end{align*}
    Which together with \zcref{eqn:BesselRootInterlacing} implies for \(x > 0\)
    \begin{equation*}
    - \infty < \difcp{\phi(\nu, x)}{x} =  \frac{W(\Jv(x), x \Jv'(x))}{(\Jv(x)^2 + (x \Jv'(x))^2)} < 0 \,.
    \end{equation*}
    Hence, \(\phi(\nu, x)\) is monotonically decreasing in \(x\). Lastly given \(x > 0, x \neq j_{\nu, k}\) and letting \(j_{\nu, 0} = 0\) then \(j_{\nu, k - 1} < x < j_{\nu, k}\) for some \(k \geq 1\) and so
    \begin{align*}
    \phi(\nu, x) & = \arctan(\nu) + \sum_{\ell = 1}^{k} \eval{\fcall{\arctan}{\frac{y \Jv'(y)}{\Jv(y)}}}_{j_{\nu, \ell - 1}+}^{\min(x, j_{\nu, \ell})-} \\
    & = \fcall{\arctan}{\frac{x \Jv'(x)}{J_\nu(x)}} - \#\{j_{\nu, k} < x\} \pi \, .
    \end{align*}
    \end{proof}
\CorBesselRootStructure
\begin{proof}
    This follows directly from \zcref{prop:BesselConditionBehaviour}.
\end{proof}

By \zcref{eqn:QuBestPossible}, \(j_{\nu, 1} \sim \nu - a_1 \spar*{\half[\nu]}^{\third}\) and so the first root in \((0,  j_{\nu, 1})\) can occur, depending on \(f\), outside the region \(0 < \nu < x\), where the uniform asymptotics \zcref{eqn:HankelUniformAsymptotic} are valid. For the problem at hand this turns out to be the case.
\PropBesselFirstRoot
\begin{proof}
    The number and conditions on the roots follows immediately from \zcref{prop:BesselConditionBehaviour}. Consider an asymptotic root, as \(\nu \to \infty\), to 
    \begin{equation*}
        \frac{v z \Jv'(v z)}{\Jv(v z)} = \nu z + \beta + \bigO(\nu z e^{- \frac{L}{R} \nu z}) \, .
    \end{equation*}
    Let \(z_\pm \coloneq {\quarter} \spar*{\sqrt{7} \pm 1}\) and \(Z \coloneq [z_-, z_+]\). Using \zcref{eqn:BesselDebyeAsymptotic} there exists a constant \(C\) such that, for \(\nu > C\) and \(z \in Z\),
    \begin{align*}
        & \Mod*{\frac{v z \Jv'(v z)}{\Jv(v z)} - \nu \spar*{1 - z^2}^{\thalf}} < C \\
        \implies & \begin{dcases}
            \frac{v z_- \Jv'(v z_-)}{\Jv(v z_-)} - \nu z_- < -\half[\nu] + C \, , \\
            \frac{v z_+ \Jv'(v z_+)}{\Jv(v z_+)} - \nu z_+ > \half[\nu] - C \, .
        \end{dcases} 
    \end{align*}
    Hence, for \(\nu > 2 C\), the first root is \(z \nu\), for some \(z \in Z\). Let \(z \sim \sum_{k = 0}^\infty c_k \nu^{-k}\). From \zcref{eqn:BesselDebyeAsymptotic}
    \begin{equation*}
        \frac{v z \Jv'(v z)}{\Jv(v z)} = v (1 - z^2)^{\frac{1}{2}} \times
    \frac{
            1 + \frac{1}{24 v}\left(3 (1 - z^2)^{-\frac{1}{2}} - 5 (1 - z^2)^{-\frac{3}{2}}\right) + \bigO(v^{-2})}{
            1 + \frac{1}{24 v}\left(-9 (1 - z^2)^{-\frac{1}{2}} + 7 (1 - z^2)^{-\frac{3}{2}}\right) + \bigO(v^{-2})
            } \,.
    \end{equation*}
    Hence,
    \begin{align*}
        & v \spar*{\sqrt{1 - c_0^2} + \bigO(v^{-1})}(1 + \bigO(v^{-1})) = \spar*{\nu (c_0 + \bigO(v^{-1})) + \bigO(1)}(1 + \bigO(v^{-1})) \\
        \implies & c_0 = \frac{1}{\sqrt{2}} \\
        \implies & \begin{aligned}[t]
            & v (\frac{1}{\sqrt{2}} + \frac{c_1}{2 v} + \bigO(v^{-2}))(1-\frac{7}{12 \sqrt{2} v} + \bigO(v^{-2})) \\
            & = (v (\frac{1}{\sqrt{2}} - \frac{c_1}{2 v}) + \half[n - 3] + \bigO(e^{-\gamma v z}))(1+\frac{5}{12 \sqrt{2} v} + \bigO(v^{-2}))
        \end{aligned} \\
        \implies & c_1 = \half[n - 2] \,.
    \end{align*}
\end{proof}

The remaining roots can be analysed with an arbitrary continuous monotonically increasing \(f\).
\PropPostTurningRoots
\begin{proof} 
    Recall from the definition of \(\theta\), \zcref{eqn:def_bessel_phase}, and \(m^2\), \zcref{eqn:def_bessel_modulus}, that
\begin{align*}
    \Jv(x) &= \sqrt{\frac{2}{\pi x}} m(\nu, x) \cos(\theta(\nu, x)), \\
    \difcp{\theta(\nu, x)}{x} &= m(\nu, x)^{-2} \, .
\end{align*}
Hence,
\begin{align*}
    \frac{x \Jv'(x)}{J_\nu(x)} & = -\half + x \frac{\difcp{m}{x}}{m} - x (\difcp{\theta}{x}) \tan(\theta) \\
    & = \frac{1}{m^2} \spar*{-\half[m^2] + \half[x] \difcp{m^2}{x} - x \tan(\theta)} \, .
\end{align*}
From \zcref{corollary:BesselRootStructure} the \(k\)-th root past \(j_{\nu, 1}\) occurs in \(\theta \in (k \pi - \half[\pi], k \pi + \half[\pi])\), and so 
\begin{align*}
    \tan(\theta - k \pi) & = \frac{1}{2}\difcp{m^2}{x} - m^2 \frac{f + \frac{1}{2}}{x} \\
    \implies k & =  
    \frac{
        \theta(\nu, x) + \arctan\left(m^2(\nu, x) \frac{f(x) + \frac{1}{2}}{x} - \frac{1}{2}\difcp{m^2(\nu, x)}{x}\right)}{\pi} \, .
\end{align*}
\end{proof}

Using \zcref{prop:PostTurningRoots} and \zcref{prop:BesselFirstRoot} the sums over \(\Sigma_{1, k}\) and \(\Sigma_{2, k}\) can be formulated as the count of positive lattice points in a function's subgraph.
\begin{proof}[Proof of \zcref{prop:bessel-sum-asymptotics}]
    For reference recall the following definitions 
    \begin{equation*}
        \xtDef \, , \xcDef \, .
    \end{equation*}
    Consider \(f(x) = \beta + x \tanh(x \oL / \R)\). First, note that 
    \begin{equation*}
        \fcall{f}{\frac{\R}{\oL} \xt(\sigma \oL)} = \sigma \R + \beta\, .
    \end{equation*}
    Using \zcref{prop:BesselFirstRoot} and \zcref{prop:PostTurningRoots}, for a fixed \(k\)
    \begin{align*}
        & \#\Set{\sigma_{1, k, \alpha} \in \Sigma_{1, k} \given \sigma_{1, k, \alpha} < \sigma} \\
        &\quad \begin{aligned}[t]
            & = \#\Set*{\alpha \given \alpha \tanh(\alpha \oL) < \sigma, \alpha \R \frac{\Jv_{k + \beta}'(\alpha \R)}{\Jv_{k + \beta}(\alpha \R)} - \beta = \alpha \R \tanh(\alpha \oL)} \\
            & = \#\Set*{\alpha \given \alpha \R  < \frac{R}{L} \xt(\sigma \oL), \alpha \R \frac{\Jv_{k + \beta}'(\alpha \R)}{\Jv_{k + \beta}(\alpha \R)} = f(\alpha \R)} \\
            & = \begin{aligned}[t]
                & \#\Set*{\alpha \given 0 < \alpha \R < \fcall{\min}{j_{k + \beta, 1}, \frac{\R}{\oL} \xt(\sigma \oL)}, \alpha \R \frac{\Jv_{k + \beta}'(\alpha \R)}{\Jv_{k + \beta}(\alpha \R)} = f(\alpha \R)} \\
                & {} + \Floor*{\frac{\fcall{\theta}{k + \beta, \frac{\R}{\oL} \xt(\sigma \oL)} + \fcall{\psi}{k + \beta, \frac{\R}{\oL} \xt(\sigma \oL), \sigma \R + \beta}}{\pi}} \, .
            \end{aligned}
        \end{aligned}
    \end{align*}
    Summing over \(k\), let
    \begin{align*}
        & S_1 \coloneq \sum_{k = 0}^\infty \omega_{n - 2, k} \#\Set*{\alpha \given 0 < \alpha \R < \fcall{\min}{j_{k + \beta, 1}, \frac{\R}{\oL} \xt(\sigma \oL)}, \alpha \R \frac{\Jv_{k + \beta}'(\alpha \R)}{\Jv_{k + \beta}(\alpha \R)} = f(\alpha \R)}\, , \\
        & S_2 \coloneq \sum_{k = 0}^\infty \omega_{n - 2, k} \Floor*{\frac{\fcall{\theta}{k + \beta, \frac{\R}{\oL} \xt(\sigma \oL)} + \fcall{\psi}{k + \beta, \frac{\R}{\oL} \xt(\sigma \oL), \sigma \R + \beta}}{\pi}} \, .
    \end{align*}
    The summands in \(S_1\) are \(1\) till \(k\) is large enough then \(0\) afterwards. Note that as \(\sigma \to \infty\)  
    \begin{equation*}
        \frac{\R}{\oL} \xt(\sigma \oL) = \sigma \R + \bigO(e^{-2 \sigma \oL}) \, ,
    \end{equation*}
    and so, from \zcref{prop:BesselFirstRoot} the summands are \(0\) for
    \begin{equation*}
        k \geq \sqrt{2} \spar*{\sigma \R + \frac{n - 2}{4}} - \beta + \bigO\spar*{\sigma^{-1}} \, .
    \end{equation*}
    Hence, for \(\sigma\) sufficiently large there exists \(C > 0\) such that 
    \begin{equation*}
        \Mod*{S_1 - \sum_{k = 0}^{\Floor*{\sqrt{2}\sigma R}} \omega_{n - 2, k}} < C \, .
    \end{equation*}
    Similarly, using \zcref{eqn:QuBestPossible}, the summands in \(S_2\) are \(0\) for
    \begin{equation*}
        k \geq \sigma R  + a_1 \spar*{\frac{\sigma R}{2}}^{\third} - \beta + \bigO(\sigma^{-\third}) \, ,
    \end{equation*}
    and for \(\sigma\) sufficiently large there exists \(C > 0\) such that
    \begin{equation*}
        \Mod*{S_2 - \sum_{k = 0}^{\Floor*{\sigma \R + a_1 \spar*{\frac{\sigma}{2}}^{\third} - \beta}} \omega_{n - 2, k} \Floor*{\frac{\fcall{\theta}{k + \beta, \frac{\R}{\oL} \xt(\sigma \oL)} + \fcall{\psi}{k + \beta, \frac{\R}{\oL} \xt(\sigma \oL), \sigma \R + \beta}}{\pi}}} < C \sigma^{-\third} \, .
    \end{equation*}
    Hence, as \(\sigma \to \infty\),
    \begin{align*}
        & \sum_{k = 0}^{\infty} \omega_{n - 2, k} \#\Set{\sigma_{1, k, \ell} \in \Sigma_{1, k} \given \sigma_{1, k, \ell} < \sigma} \\
        & = S_1 + S_2 \\
        & = \begin{aligned}[t]
            \sum_{k = 0}^{\Floor*{\sigma \R + a_1 \spar*{\frac{\sigma}{2}}^{\third} - \beta}} \omega_{n - 2, k} \Floor*{\frac{\fcall{\theta}{k + \beta, \frac{\R}{\oL} \xt(\sigma \oL)} + \fcall{\psi}{k + \beta, \frac{\R}{\oL} \xt(\sigma \oL), \sigma \R + \beta}}{\pi}} \\
            {} + \sum_{k = 0}^{\Floor*{\sqrt{2}\sigma R}} \omega_{n - 2, k} + \bigO(1) \, .
        \end{aligned}\
    \end{align*}
    
    Finally, comparing \(\Sigma_{1, k}\) to \(\Sigma_{2, k}\), the previous analysis carries over mutatis mutandis with \(\xt\) replaced with \(\xc\). There are a finite number of first roots not present in \(\Sigma_{2, k}\), explicitly when \(k \leq \frac{\R}{\oL}\) and as \(\sigma \to \infty\)
    \begin{equation*}
        \frac{R}{L} \xc(\sigma \oL) = \frac{R}{L} \xt(\sigma \oL) + \bigO(e^{-2 \sigma \oL}) \, ,
    \end{equation*}
    hence
    \begin{equation*}
        \sum_{k = 0}^{\infty} \omega_{n - 2, k} \#\Set{\sigma_{1, k, \ell} \in \Sigma_{1, k} \given \sigma_{1, k, \ell} < \sigma} = \sum_{k = 0}^{\infty} \omega_{n - 2, k} \#\Set{\sigma_{2, k, \ell} \in \Sigma_{2, k} \given \sigma_{2, k, \ell} < \sigma} + \bigO(1) \, .
    \end{equation*}
\end{proof}

The summand, \(\theta + \psi\), is not easily summable, and so a more elementary approximation is sought. First the asymptotic behaviour of the constituent parts, \(\theta\), \(m^2\) and \(\difcp{m^2}{x}\), are analysed for the region \(0 < \nu < x\) as \(x \to \infty\).
\PropBesselCountingApproximations
\begin{proof}
    Recall from \zcref{eqn:BesselPhaseBounds}
    \begin{align*}
        0 \leq \tilde{\theta}(\nu, x) - \theta(\nu, x) \leq \BesselPhaseApproxErrorDef \, .
    \end{align*}
    Hence, the first inequality follows from noting that
    \begin{align*}
        \frac{3 x^2 + 2 \nu^2}{24 (x^2 - \nu^2)^{\half[3]}} = (x - \nu)^{-\half[3]}\frac{3 x^2 + 2 \nu^2}{24 (x + \nu)^{\half[3]}} \leq \frac{1}{8}\spar*{\frac{x^{\third}}{x - v}}^{\half[3]}\, .
    \end{align*}
    Next, using the uniform Hankel expansions \zcref{eqn:HankelUniformAsymptotic} there exists a constant \(C > 0\) such that, for \(x > C\) and \(0 < \nu < x\),
    \begin{equation*}
        \Mod*{\frac{\Hv^{j}(x) \nu^{\third}}{\spar*{\frac{4 \xi \nu^2}{x^2 - \nu^2}}^{\quarter} \spar*{
            \fcall{\Ai}{- \nu^{\third[2]} \xi} \pm \iu \fcall{\Bi}{- \nu^{\third[2]} \xi}}} - 1} < C \nu^{-1}\min(1, \xi^{-\half[3]}) \, .
    \end{equation*}
    Hence, with \(M_A\) as in \zcref{eqn:AiryModulusAsymptotic}, there exists a constant \(C'\) such that, for \(x > C' > C\) and \(0 < \nu < x\),
    \begin{equation*}
        \Mod*{\frac{\Hv^{1}(x) \Hv^{2}(x) \nu^{\third[2]}}{\spar*{\frac{4 \xi \nu^2}{x^2 - \nu^2}}^{\half}M_A(-\nu^{\third[2]}\xi)^2} - 1} < 2 C \nu^{-1}\min(1, \xi^{-\half[3]}) + C^2 \nu^{-2}\min(1, \xi^{-3}) < C' \nu^{-1}\min(1, \xi^{-\half[3]}) \, .
    \end{equation*}
    The second inequality now follows from
    \begin{align*}
        & m^2(\nu, x) = \frac{\pi x}{2} H^{(1)}_\nu(x) H^{(2)}_\nu(x) \, , \\
        & \pi \nu^{\third} \xi^{\half} M_A(- \nu^{\third[2]} \xi)^2 = m(\tthird, \nu \xi^{\half[3]})^2 \, .
    \end{align*}
    The third inequality follows similarly to the second, using \zcref{eqn:HankelUniformAsymptotic}, \zcref{eqn:DHankelUniformAsymptotic}, and
    \begin{align*}
        \difcp{m^2}{x}
            & = \frac{m^2}{x} + \frac{\pi x}{2} \Re(H_\nu^{(1)} \difcp{H_\nu^{(2)}}{x}) \\
            & = \frac{m^2}{x} + \frac{\pi}{2} \spar*{\diff*{M_A^2}{x}}(-\nu^{\third[2]}\xi)\spar*{1 + \bigO(\nu^{-1}\min(1, \xi^{-\half[3]}))} \, .
    \end{align*}
\end{proof}
\begin{proof}[Proof of \zcref{prop:bessel-lattice-point-func-approx}]
    For reference recall the following definition
    \begin{equation*}
        \EtaOneDef \, ,
    \end{equation*}
    together with the claim to be proved that given \(\third < \gamma \leq 1\) then, uniformly in \(\beta < k + \beta < \sigma \R - {\thalf}(\sigma \R)^\gamma\), as \(\sigma \to \infty\),
    \begin{equation*}
        \fcall{\theta}{k + \beta, \frac{R}{L} \xt(\sigma \oL)} + \fcall{\psi}{k + \beta, \frac{R}{L} \xt(\sigma \oL), \sigma \R + \beta} = \eta_1(k + \beta, \sigma \R) - \quarter[\pi] + \fcall{\bigO}{\sigma^{\half[1 - 3 \gamma]}} \, .
    \end{equation*}
    First, the elementary approximations for \(\theta\), \(m^2\) and \(\difcp{m^2}{x}\) are \(\tilde{\theta}\), \(\spar*{1 - \spar*{\frac{\nu}{x}}^2}^{-\half}\) and \(0\) respectively. However, these approximations start breaking down when \(x - v \approx x^{\third}\). From \zcref{prop:BesselCountingApproximations} there exists a constant \(C > 0\) such that uniform limits in terms of \(C\) apply for \(x > C\). Hence, consider the parameterised regions for \({\tthird} < \gamma \leq 1\)
    \begin{equation*}
      \RIDef \,.
    \end{equation*} The following bounds apply in \(R(\gamma, C)\):
    \begin{align*}
        & 0 < \varepsilon_\theta < \frac{1}{8} x^{\half[1 - 3 \gamma]} \, , & \text{\zcref{prop:BesselCountingApproximations}}\\
        & \half[1] x^{\half[3 \gamma - 1]} < \nu \xi^{\half[3]} < \half[3] x \, , & \text{\zcref{eqn:xiDef}}\\
        & 1 < \spar*{1 - \spar*{\frac{\nu}{x}}}^{-\half} < \frac{2}{\sqrt{3}} x^{\half[1 - \gamma]} \, , \\
        & 0 < 1 - m({\tthird}, \nu \xi^{\half[3]})^2 < \frac{5}{8} x^{1 - 3 \gamma} \, , & \text{\zcref{eqn:AiryModulusAsymptotic}}\\
        & \Mod*{m(\nu, x)^2 \spar*{1 - \spar*{\frac{\nu}{x}}}^{\half} - 1} < \frac{5}{8} x^{1 - 3 \gamma} + 2 C x^{-1} \, , & \text{\zcref{prop:BesselCountingApproximations}}\\
        & \frac{1}{3 \pi} x^{-1} - \frac{35}{216} x^{-3}  < \spar*{M_A^2}'(- \nu^{\third[2]} \xi) < \frac{1}{\pi} x^{\half[1 - 3 \gamma]}   \, . & \text{\zcref{eqn:DAiryModulusAysmptotic}}
    \end{align*}
    Next recall that
    \begin{align*}
        - \half[\pi] < \arctan(u) - \arctan(v) < \half[\pi] & \iff - \half[\pi] < \arctan(u^{-1}) - \arctan(v^{-1}) < \half[\pi] \\ 
        & \implies \arctan(u) - \arctan(v) = \fcall{\arctan}{\frac{u - v}{1 + u v}}\, .
    \end{align*}
    For arbitrary \(c_0, c_1 \in \mathbb{R}\) consider
    \begin{align*}
        & u = m(\nu, x)^2 + c_0 x^{-1} m(\nu, x)^2 + c_1 \spar*{\difcp{m(\nu, x)^2}{x} - x^{-1}m(\nu, x)^2} \, , \\
        & v = \spar*{1 - \spar*{\frac{\nu}{x}}^2}^{-\half}\, .
    \end{align*}
    From the above bounds, \zcref{prop:BesselCountingApproximations}, and with abuse of notation, there exists a constant \(C > 0\) depending on \(\gamma, c_0, c_1\) such that in \(R(\gamma, C)\)
    \begin{equation*}
        \Mod*{\frac{u}{v} - 1} < C x^{\half[1 - 3 \gamma]} < {\half} \implies - \frac{\pi}{2} < \arctan(u) - \arctan(v) < \frac{\pi}{2} \,.
    \end{equation*}
    Moreover,
    \begin{align*}
        & \Mod*{\frac{u - v}{1 + u v}} = v^{-1} \Mod*{\frac{\frac{u}{v} - 1}{\frac{u}{v} - 1 + 1 + v^{-2}}} \leq 2 v^{-1} \Mod*{\frac{u}{v} - 1} < 2 C x^{\half[1 - 3 \gamma]} \, , \\
        \implies & \Mod*{\arctan(u) - \arctan(v)} <  2 C x^{\half[1 - 3 \gamma]} \, .
    \end{align*}    
    Lastly, let \(x = \sigma R, \nu = k + \beta, c_0 = \beta + {\half}, c_1 = {\half}\) and recall that \(\frac{\R}{\oL}\xt(\sigma \oL) = \sigma \R + \bigO(e^{-2 \sigma L})\). Hence,
    \begin{align*}
        & \fcall{\theta}{k + \beta, \frac{R}{L} \xt(\sigma \oL)} + \fcall{\psi}{k + \beta, \frac{R}{L} \xt(\sigma \oL), \sigma \R + \beta} \\
        & \quad \begin{aligned}
            & = \tilde{\theta}(k + \beta, \sigma \R) + \fcall{\arctan}{\frac{\sigma R}{\sqrt{(\sigma \R)^2 - (k + \beta)^2}}} + \fcall{\bigO}{\sigma^{\half[1 - 3 \gamma]}} \\
            & = \eta_1(k + \beta, \sigma \R) - \quarter[\pi] + \fcall{\bigO}{\sigma^{\half[1 - 3 \gamma]}} \, ,
        \end{aligned}
    \end{align*}
    as required
\end{proof}

\begin{section}{Lattice Point Count Formulation Over Radially Localised Waves}\label{section:lattice-point-mod-bessel}
The next step is to consider sums containing \(\Sigma_{3, k}\) and \(\Sigma_{4, k}\). These correspond to radially localised waves, in that the eigenfunctions are of the form \(u(r, x, z) = I_\nu(\alpha r) \dots\), and so, exponentially decay as \(r \to 0\). From \zcref{prop:problem-solns} the associated transcendental equations are of the form
\begin{align*}
    & \alpha \frac{\Iv_{k + \beta}'(\alpha R)}{\Iv_{k + \beta}(\alpha R)} - \frac{\beta}{R} = - \alpha \tan(\alpha L) \, ,\\
    & \alpha \frac{\Iv_{k + \beta}'(\alpha R)}{\Iv_{k + \beta}(\alpha R)} - \frac{\beta}{R} = \alpha \cot(\alpha L) \, .
\end{align*}
Again, the right-hand side is a piecewise smooth function with multi-valued inverse and the left-hand side is a monotonically increasing smooth function. However, in contrast to transversely localised waves, both \(\tan\) and \(\cot\) have a trivial phase formulation and the complexity is in the analysis of the ratio of modified Bessel functions. First, counting the roots can be formulated as the count of positive lattice points in a function's subgraph:
\PropModifiedBesselCount
\begin{proof}
    Consider the possible roots for the equation
    \begin{equation*}
        \frac{x \Iv'(x)}{\Iv(x)} + c = -x \tan(\mu x + \vartheta) \, ,
    \end{equation*}
    with \(\vartheta \in [0, \half[\pi]]\). Note that the right-hand side is monotonically decreasing between its poles, starting from
    \begin{align*}
        \lim_{x \to 0+} -x \tan(\mu x + \vartheta) = \frac{1}{\mu} \mathbb{1}_{\Set*{\vartheta = \thalf[\pi]}} \, .
    \end{align*}
    Whilst the left-hand side is monotonically increasing from
    \begin{align*}
        \lim_{x \to 0+} \frac{x \Iv'(x)}{\Iv(x)} + c = \nu + c \, .
    \end{align*}
    First consider the range \(\vartheta \in [0, \thalf[\pi])\), there are no roots in \((0, \half[\pi] - \vartheta]\), one possible root at \(0\), if \(\nu + c = 0\), and a single root in the range \((k \pi - \half[\pi] - \vartheta, k \pi + \half[\pi] - \vartheta)\), for each \(k \in \mathbb{N}\). So at a non-zero root
    \begin{align*}
        & \frac{x \Iv'(x)}{\Iv(x)} + c = -x \tan(\mu x + \vartheta + k \pi) \\
        & \implies k = \frac{\mu x + \vartheta + \arctan\left(\frac{\Iv'(x)}{\Iv(x)} + \frac{c}{x}\right)}{\pi} \, .
    \end{align*}
    Furthermore,
    \begin{equation*}
        \frac{\Iv'(x)}{\Iv(x)} + \frac{c}{x} = \frac{\nu + c}{x} + \bigO(x) \implies \lim_{x \to 0+} \arctan\left(\frac{\Iv'(x)}{\Iv(x)} + \frac{c}{x}\right) = \thalf[\pi]\mathbb{1}_{\Set*{\nu + c \neq 0}} \, ,
    \end{equation*}
    and so \(k\) corresponds to the \(k\)-th non-zero root.
    
    The analysis is similar for \(\vartheta = \half[\pi]\), however there are additional roots for the range \(0 \leq \nu + c < 1\) before the initial root becomes a root at zero, and so, for \(\nu + c \geq 1\), \(k\) corresponds to the \(k - 1\)-th non-zero root.
\end{proof}
\begin{proof}[Proof of \zcref{prop:mod-bessel-sum-asymptotics}]
    For reference, recall the following definition \(\xvDef\). Using \zcref{prop:ModifiedBesselCount} for a fixed \(k\) with \(x = \alpha R\), \(y = \xv_{k + \beta}((\sigma + \beta) R)\) and \(\mu = \frac{\oL}{\R}\)
    \begin{align*}
        & \#\Set{\sigma_{3, k, \ell} \in \Sigma_{3, k} \given \sigma_{3, k, \ell} < \sigma} \\
        &\quad \begin{aligned}[t]
            & = \#\Set*{\alpha \given \alpha \frac{\Iv_{k + \beta}'(\alpha R)}{\Iv_{k + \beta}(\alpha R)} - \beta < \sigma, \alpha \R \frac{\Iv_{k + \beta}'(\alpha R)}{\Iv_{k + \beta}(\alpha R)} - \beta = - \alpha \R \tan(\alpha L)} \\
            & = \#\Set*{\alpha \given \alpha \R  < \xv_{k + \beta}(\sigma \R + \beta), \alpha \R \frac{\Iv_{k + \beta}'(\alpha R)}{\Iv_{k + \beta}(\alpha R)} - \beta = - \alpha \R \fcall{\tan}{\frac{\oL}{\R} \alpha \R}} \\
            & = \Floor*{\frac{\frac{L}{R} \xv_{k + \beta}(\sigma \R + \beta) + \fcall{\arctan}{\frac{\sigma \R}{\xv_{k + \beta}(\sigma \R + \beta)}}}{\pi}} \, .
        \end{aligned}
    \end{align*}
    From \zcref{eqn:MonotonicityModifiedConditionFunc}
    \begin{equation*}
        \xv_{\sigma \R + \beta}(\sigma \R + \beta) = 0 \, ,
    \end{equation*}
    and so 
    \begin{align*}
        & \sum_{k = 0}^{\infty} \omega_{n - 2, k} \#\Set{\sigma_{3, k, \alpha} \in \sigma_{3, k} \given \sigma_{3, k, \alpha} < \sigma} \\
        & = \sum_{k = 0}^{\Floor*{\sigma R}} \omega_{n - 2, k} \Floor*{\frac{\frac{L}{R} \xv_{k + \beta}(\sigma R + \beta) + \fcall{\arctan}{\frac{\sigma R}{\xv_{k + \beta}(\sigma R + \beta)}}}{\pi}} \, ,
    \end{align*}
    as required. The sum over \(\Sigma_{4, k}\) is shown similarly.
\end{proof}

Last, the asymptotics of \(\xv_\nu\), which lead to an approximation of the summand.
\PropModifiedBesselInverseAsymptotics
\begin{proof}
    \(X_\nu\) will be approximated via a sequence of uniform approximations. Each approximation will be stated as \(f(x, \nu, \sigma) = \bigO(g(\sigma))\), meaning that there exist constants \(C_0\) and \(C_1\) such that, for \(\sigma > C_0\) and \(\sqrt{x^2 + \nu^2} \geq \sigma\) then \(\Mod*{f(x, \nu, \sigma)} < C_1 g(\sigma)\). The approximation of \(\frac{x I'_\nu(x)}{I_\nu(x)}\) follow from combining \zcref{eqn:ModifiedBesselAsymptotics}, \zcref{eqn:DModifiedBesselAsymptotics}, in particular
    \begin{equation*}
        \Mod*{\frac{x I'_\nu(x)}{I_\nu(x)} \frac{p(x, \nu)}{\nu} - \frac{1 + \frac{1}{24 \nu}\spar*{-9 p(x, \nu) + 7 p(x, \nu)^3}}{1 + \frac{1}{24 \nu}\spar*{3 p(x, \nu) - 5 p(x, \nu)^3}}} = \bigO\spar*{\spar*{x^2 + \nu^2}^{-1}} \, .
    \end{equation*}
    By \zcref{eqn:ModifiedConditionFuncBound} \(X_\nu(\sigma)\) is within the region of validity, \(\sigma \leq \sqrt{x^2 + \nu^2} < \infty\), for the uniform asymptotics. It is convenient to re-parameterise the approximation, let \(s(\nu, \sigma) \coloneq \frac{\nu}{\sigma}\) and \(y(\nu, \sigma) \coloneq \frac{p(X_\nu(\sigma), \nu)}{s}\),
    \begin{align*}
        & \sigma \leq \sqrt{X_\nu(\sigma)^2 + \nu^2} < \infty \iff 0 < p \leq s \iff 0 < y \leq 1, 0 < s \leq 1 \, , \\
        & X_\nu(\sigma) = \sigma \sqrt{y^{-2} - s^2} \, , \\
        & \Mod*{y - \frac{1 + \frac{1}{24 \sigma}\spar*{-9 y + 7 s^2 y^3}}{1 + \frac{1}{24 \sigma}\spar*{3 y - 5 s^2 y^3}}} = \bigO\spar*{\sigma^{-2} y^{2}} \, .
    \end{align*}
    Using standard root counting algorithms the polynomial
    \begin{equation*}
        y \spar*{1 + \frac{1}{24 \sigma}\spar*{3 y - 5 s^2 y^3}} - \spar*{1 + \frac{1}{24 \sigma}\spar*{-9 y + 7 s^2 y^3}}
    \end{equation*}
    has a single root in \(0 \leq y \leq 1\), denoted \(\tilde{y}\). Furthermore, for \(\sigma > \frac{1}{12}\), \(0 < y \leq 1\) and \(0 < s \leq 1\)
    \begin{align*}
        & 1 - \frac{1}{2 \sigma} < \frac{1 + \frac{1}{24 \sigma}\spar*{-9 y + 7 s^2 y^3}}{1 + \frac{1}{24 \sigma}\spar*{3 y - 5 s^2 y^3}} < 1 \, , \\
        & \Mod*{\difcp{\frac{1 + \frac{1}{24 \sigma}\spar*{-9 y + 7 s^2 y^3}}{1 + \frac{1}{24 \sigma}\spar*{3 y - 5 s^2 y^3}}}{y}} \leq \frac{12}{12 \sigma - 1} \, .
    \end{align*}
    Hence, \(y\) can be approximated by \(\tilde{y}\),
    \begin{align*}
        & \Mod*{y - \tilde{y}} \leq \Mod*{y - \frac{1 + \frac{1}{24 \sigma}\spar*{-9 y + 7 s^2 y^3}}{1 + \frac{1}{24 \sigma}\spar*{3 y - 5 s^2 y^3}}} + \frac{12}{12 \sigma - 1} \Mod*{y - \tilde{y}} \\
        \implies & \Mod*{y - \tilde{y}} = \bigO\spar*{\sigma^{-2}} \\
        \implies & \Mod*{X_\nu(\sigma) - \sigma \sqrt{\tilde{y}^{-2} - s^2}} = \Mod*{\sigma \sqrt{y^{-2} - s^2} - \sigma \sqrt{\tilde{y}^{-2} - s^2}} = \bigO(\sigma^{-1}) \, .
    \end{align*}
    Next, \(\tilde{y}\) can be approximated by \(1 - \frac{1 - s^2}{2 \sigma}\),
    \begin{align*}
        & \Mod*{\tilde{y} - \spar*{1 - \frac{1 - s^2}{2 \sigma}}} = \bigO(\sigma^{-2}) \\
        \implies & \Mod*{\sigma \sqrt{\tilde{y}^{-2} - s^2} - \sigma \spar*{{\spar*{1 - \frac{1 - s^2}{2 \sigma}}^{-2} - s^2}}^{\half}} = \bigO(\sigma^{-1}) \, .
    \end{align*}
    Finally,
    \begin{equation*}
        \Mod*{\sigma \spar*{{\spar*{1 - \frac{1 - s^2}{2 \sigma}}^{-2} - s^2}}^{\half} - \sqrt{1 - s^2}(\sigma + {\thalf})} = \bigO(\sigma^{-1}) \, .
    \end{equation*}
    Hence,
    \begin{equation*}
        \Mod*{X_\nu(\sigma) - \sqrt{1 - s^2}(\sigma + {\thalf})} = \Mod*{X_\nu(\sigma) - \sqrt{\sigma^2 - \nu^2}\spar*{1 + \frac{1}{2 \sigma}}} = \bigO(\sigma^{-1}) \, .
    \end{equation*}
\end{proof}
\begin{proof}[Proof of \zcref{prop:mod-bessel-lattice-point-func-approx}]
    For reference recall the definition
    \begin{equation*}
        \EtaTwoDef \, .
    \end{equation*}
    Using \zcref{prop:ModifiedBesselInverseAsymptotics}, then uniformly in \(0 < \nu < \sigma\), as \(\sigma \to \infty\),
    \begin{align*}
        & \frac{L}{R}\xv_{\nu}(\sigma) + \fcall{\arctan}{\frac{\sigma - \beta}{\xv_{\nu}(\sigma)}} \\
        & \qquad \begin{aligned}
        & = \sqrt{\sigma^2 - \nu^2}\spar*{1 + \frac{1}{2 \sigma}} + \fcall{\arctan}{\frac{\sigma - \beta}{\sqrt{\sigma^2 - \nu^2}}\spar*{1 + \frac{1}{2 \sigma}}} + \bigO(\sigma^{-1}) \\
        & = \sqrt{\sigma^2 - \nu^2}\spar*{1 + \frac{1}{2 \sigma}} + \half[\pi] - \fcall{\arctan}{\frac{\sqrt{\sigma^2 - \nu^2}}{\sigma}} + \bigO(\sigma^{-1}) \\
        & = \eta_2(\nu, \sigma) + \half[\pi] + \bigO(\sigma^{-1})\, .
        \end{aligned}
    \end{align*}
\end{proof}

\end{section}

\section{Lattice Point Count Approximation}\label{section:lattice-point-count-approx}

The last step is to approximate the lattice point counts. To prove \zcref{prop:BesselStrongConvexity}, \zcref{prop:ModifiedBesselStrongConcavity} and \zcref{thm:main} three points need to be addressed: 
\begin{itemize}
    \item the breakdown of the elementary approximations \(\eta_1(k + \beta, \sigma R)\) when \(k + \beta \approx \sigma R\),
    \item bounding \(\sum_{k = 0}^{\ell} \omega_{n - 2, k} \spar*{\Floor{\eta_j(k) + c} - \eta_j(k) - c}\),
    \item calculating \(\sum_{j = 0}^{\ell} \omega_{n - 2, k} \spar*{\eta_j(k) + c}\),
\end{itemize}
considered in that order.

First, the breakdown of the elementary approximations can be dealt with via truncation.
\PropTruncateBesselSum
\begin{proof}
    Denote \(\sigma_- = \sigma \R - {\thalf} (\sigma \R)^{\gamma} + \beta\) and \(\sigma_+ = \sigma \R - {\thalf} (\sigma \R)^{\third} + \beta\). 
    By definition
    \begin{align*}
        & -\half[\pi] \leq \psi(k + \beta, \frac{\R}{\oL} \xt(\sigma \oL), \sigma \R + \beta) \leq \half[\pi] \, ,\\
        & 0 < \omega_{n - 2, k} \leq \omega_{n - 2, \Floor*{\sigma_+}} \, .
    \end{align*}
    Next, from \cite[\S 13.73 eq. (2)]{WatsonBessel},
    \begin{align*}
        \difcp{\theta(\nu, x)}{\nu} = - \frac{4}{\pi M^2(\nu, x)} \int_{0}^{\infty} K_0(2 x \sinh(t)) e^{-2 \nu t} \dl{t} < 0 \, ,
    \end{align*}
    and so for \(\sigma_- < k + \beta \leq \sigma_+\)
    \begin{equation*}
        -\half[\pi] \leq \theta(k + \beta, \frac{\R}{\oL} \xt(\sigma \oL)) \leq \theta(\sigma_-, \frac{\R}{\oL} \xt(\sigma \oL)) \leq \widetilde{\theta}(\sigma_-, \sigma \R) \, .
    \end{equation*}
    Hence,
    \begin{align*}
        & \widetilde{\theta}(\sigma_-, \sigma \R) + \quarter[\pi] = \bigO\left(\left(\frac{\sigma^{1 + \gamma}}{\sigma^{\third[4]}}\frac{(2 - \sigma^{\gamma - 1})}{(1 - \sigma^{\gamma - 1})^{\third[4]}} \right)^{\half[3]}\right) = \bigO(\sigma^{\half[3 \gamma] - \half}) \, , \\
        & \sigma_- - \sigma_+ = \bigO(\sigma^\gamma) \, , \\
        & \omega_{n - 2, \Floor*{\sigma_+}} = \bigO(\sigma^{\n - 3}) \, .
    \end{align*}
    Hence,
    \begin{align*}
        & \sum_{\sigma_- < k \leq \sigma_+} \omega_{n - 2, k} \Floor*{\frac{\eta(k + \beta, \frac{\R}{\oL} \xt(\sigma \oL)) + \psi(k + \beta, \frac{\R}{\oL} \xt(\sigma \oL), (\sigma + \beta) \R)}{\pi}} \\
        & \qquad \begin{aligned}[t]
            & \leq \spar*{\frac{\widetilde{\eta}(\sigma_-, \sigma \R)}{\pi} + \thalf[3]} \spar*{\sigma_- - \sigma_+} \omega_{n - 2, \Floor*{\sigma}} \\
            & = \bigO(\sigma^{\half[5 \gamma] - \half} + \sigma^\gamma) \sigma^{\n - 3} \, .
        \end{aligned}
    \end{align*}
\end{proof}

Next, consider bounding the error associated with taking the integer part of a function. Applying \zcref{prop:VanDerCorputErrorEstimate} and \zcref{eqn:euler-maclaurin} requires checking strong convexity of \(\eta_1\), checking strong concavity of \(\eta_2\) and bounding derivatives. Both proofs use the following statement
\begin{equation*}
    \RationalFuncBounds \,.
\end{equation*}
\begin{proof}[\zcref{prop:BesselStrongConvexity} Proof]
    For reference recall the definition
    \begin{equation*}
        \omega_{1, k} \coloneq \begin{dcases}
            1 & \mid k = 0 \, , \\
            2 & \mid \mathrm{ otherwise} \, ,
        \end{dcases}
    \end{equation*}
    and the claim, to be proved, that given \(c \in \mathbb{R}\) then uniformly in \(k\) as \(\sigma \to \infty\)
    \begin{equation*}
        \sum_{k \leq j \leq \sigma} \omega_{1, j - k} \Floor*{\frac{\eta_1(k, \sigma)}{\pi} + c} = 2 \int_{k}^{\sigma} \frac{\eta_1(\nu, \sigma)}{\pi} \dl{\nu} - (1 - 2c)(\sigma - k) + \bigO(\sigma^{\third[2]}) \, .
    \end{equation*}
    The plan is to use \zcref{prop:VanDerCorputErrorEstimate} and the Euler-Maclaurin formula, which requires strong convexity of \(\eta_1\). First, for \(z \in [0, 1]\),
    \begin{equation*}
        \spar*{\difcp[2]{\eta_1}{\nu}}(z \sigma, \sigma) = \frac{1}{\sqrt{1-z^2} \sigma} + \frac{- 2 z^4 + z^2 + 2}{(2 - z^2)^2} \frac{1}{\left(1-z^2\right)^{3/2}\sigma^2} > \frac{1}{\sigma} \, .
    \end{equation*}
    Hence, \(\eta_1\) is strongly convex, however \(\difcp{\eta_1}{\nu}(\nu, \sigma)\) is singular as \(\nu \to \sigma\), so a truncated range is sought where \(\difcp{\eta_1}{\nu}\) is bounded. In particular for \(\sigma > 1\),
    \begin{align*}
        & \spar*{\difcp{\eta_1}{\nu}}(\sigma - 1, \sigma) = \frac{1}{\sqrt{2 \sigma - 1}}\spar*{1 - \frac{3 \sigma - 1}{\sigma(\sigma + 1) + (\sigma - 1)}} - \fcall{\arctan}{\frac{\sqrt{2 \sigma - 1}}{\sigma - 1}} \\
        \implies & - \half[\pi] < - \fcall{\arctan}{\frac{\sqrt{2 \sigma - 1}}{\sigma - 1}} < \spar*{\difcp{\eta_1}{\nu}}(\sigma - 1, \sigma) < \frac{1}{\sqrt{2 \sigma - 1}} < 1 \, .
    \end{align*} 
    Hence, for all \(0 < k \leq \sigma - 1\), by \zcref{prop:VanDerCorputErrorEstimate}
    \begin{equation*}
        \Mod*{\sum_{k < j \leq \sigma - 1} \rho_1(\eta_1(j, \sigma))} \leq \half[11 \pi] \sigma^{\third[2]} + 11 \sigma^{\thalf} = \fcall{\bigO}{\sigma^{\third[2]}} \, .
    \end{equation*}
    Finally, using the Euler-Maclaurin formula,
    \begin{align*}
        \sum_{k < j \leq \sigma} \Floor*{\frac{\eta_1(j, \sigma)}{\pi} + c} 
            & = \sum_{k < j \leq \sigma} \spar*{\frac{\eta_1(j, \sigma)}{\pi} + c} - \half[\Floor*{\sigma} - k] + \sum_{k < j \leq \sigma} \fcall{\rho_1}{\frac{\eta_1(j, \sigma)}{\pi} + c} \\
            & = \begin{aligned}[t]
                \int_{k}^{\sigma - 1} \frac{\eta_1(\nu, \sigma)}{\pi} + c \dl{\nu} - \half[\sigma - k] - \eval{\spar*{\rho_1(\nu)\spar*{\frac{\eta_1(\nu, \sigma)}{\pi} + c}}}_{\nu = k}^{\nu = (\sigma - 1)-} \\
                {} + \bigO(1 + \sigma^{\half[2]} + \Mod*{\eval{\difcp{\eta_1}{\nu}}_{\nu = k}^{\nu = \sigma - 1}} \sigma^{\third[2]})
            \end{aligned} \\
            & = \int_{k}^{\sigma} \frac{\eta_1(\nu, \sigma)}{\pi}\dl{\nu} - (1 - 2c)\half[\sigma - k] - \frac{\eta_1(k, \sigma)}{2 \pi} + \bigO(\sigma^{\third[2]}) \,.
    \end{align*}
    This implies the weighted sum on noting that
    \begin{align*}
        \sum_{k \leq j \leq \sigma} \omega_{1, j} \Floor*{\frac{\eta_1(j, \sigma)}{\pi} + c} = \Floor*{\frac{\eta_1(k, \sigma)}{\pi} + c} + 2 \sum_{k < j \leq \sigma} \Floor*{\frac{\eta_1(j, \sigma)}{\pi} + c} \, .
    \end{align*}
\end{proof}
\begin{proof}[\zcref{prop:ModifiedBesselStrongConcavity} Proof]
    For reference recall the claim, to be proved, that given \(c \in \mathbb{R}\) then uniformly in \(k\) as \(\sigma \to \infty\)
    \begin{equation*}
        \sum_{k \leq j \leq \sigma} \omega_{1, j - k} \Floor*{\frac{\eta_2(j, \sigma)}{\pi} + c} = 2 \int_{k}^{\sigma} \frac{\eta_2(\nu, \sigma)}{\pi} \dl{\nu} - (1 - 2c)(\sigma - k) + \bigO(\sigma^{\third[2]}) \, .
    \end{equation*}
    To apply \zcref{prop:VanDerCorputErrorEstimate} requires strong concavity of \(\eta_2\), and also splitting the range of summation such that the derivative is bounded. For fixed \(\sigma\) let \(\mu_2(\eta_2(\nu, \sigma)) = \nu\). First, for \(z \in [0, 1]\),
    \begin{align*}
        & \spar*{\difcp{\eta_2}{\nu}}(z \sigma, \sigma) = -\frac{z}{\sqrt{1 - z^2}} \spar*{ \frac{\oL}{\R} - \sigma^{-1}\spar*{\frac{1}{(2 - z^2)} -  \frac{\oL}{2 \R}}} \, ,\\
        & \spar*{\difcp[2]{\eta_2}{\nu}}(z \sigma, \sigma) = - \sigma^{-1} (1 - z^2)^{-\half[3]} \spar*{ \frac{\oL}{\R} - \sigma^{-1}\spar*{\frac{- 2 z^4 + z^2 + 2}{(2 - z^2)^2} -  \frac{\oL}{2 \R}}} \, , \\
        & \spar*{\difcp[2]{\mu_2}{\nu}}(\eta_2(z \sigma, \sigma)) = - \sigma^{-1} z^{-3} \frac{ \frac{\oL}{\R} - \sigma^{-1}\spar*{\frac{- 2 z^4 + z^2 + 2}{(2 - z^2)^2} -  \frac{\oL}{2 \R}}}{\spar*{ \frac{\oL}{\R} - \sigma^{-1}\spar*{\frac{1}{(2 - z^2)} -  \frac{\oL}{2 \R}}}^3} \, .
    \end{align*}
    Consider the range \(0 \leq z \leq \frac{1}{\sqrt{2}}\). For \(\sigma > 2 \frac{\R}{\oL}\),
    \begin{align*}
        & \spar*{\difcp[2]{\eta_2}{\nu}}(z \sigma, \sigma) < - \sigma^{-1} \frac{15}{32}  \frac{\oL}{\R} \, ,\\
        & -  \frac{\oL}{\R} - \spar*{ \frac{\oL}{2 \R}}^2 < \spar*{\difcp{\eta_2}{\nu}}(z \sigma, \sigma) \leq 0 \, .
    \end{align*}
    Next, consider the range \(\frac{1}{\sqrt{2}} \leq z \leq 1\),
    \begin{align*}
        & \spar*{\difcp[2]{\mu_2}{\nu}}(\eta_2(z \sigma, \sigma)) < - \sigma^{-1} \frac{\frac{15}{32}  \frac{\oL}{\R}}{\spar*{ \frac{\oL}{\R} + \spar*{ \frac{\oL}{2 \R}}^2}^3} \, ,\\
        & - \frac{3}{2}\frac{\R}{\oL} < \spar*{\difcp[1]{\mu_2}{\nu}}(\eta_2(z \sigma, \sigma)) \leq 0 \, .
    \end{align*}
    Hence, by \zcref{prop:VanDerCorputErrorEstimate} there exists a constant \(C\) such that, for \(\sigma > C\) and \(k > 0\), 
    \begin{align*}
        & \Mod*{\sum_{k < j \leq \frac{\sigma}{\sqrt{2}}} \rho_1(\eta_2(j, \sigma))} <  C \sigma^{\third[2]} \, , \\
        & \Mod*{\sum_{k < j \leq \fcall{\eta_2}{\frac{\sigma}{\sqrt{2}}, \sigma}} \rho_1(\mu_2(j, \sigma))} < C \sigma^{\third[2]} \, .
    \end{align*}
    For the weighted sum, to ease notation, fix \(\sigma\), let \(f(\nu) = \frac{\eta_2(\nu, \sigma)}{\pi} + c\) and \(g(f(\nu)) = \nu\). Note that
    \begin{align*}
        & \sum_{k < j \leq \sigma} \Floor*{f(j)} = \sum_{f(\sigma) < j \leq f(k)} \Floor*{g(j)} + \Floor{f(\sigma)}\spar*{\Floor{\sigma} - k} - k \spar*{\Floor{f(k)} - \Floor{f(\sigma)}} \, , \\
        & \int_{k}^{\sigma} f(\nu) \dl{\nu} = \int_{f(\sigma)}^{f(k)} g(\omega) \dl{\omega} + f(\sigma)(\sigma - k) - k(f(k) - f(\sigma)) \, .
    \end{align*}
    Hence, if \(k \geq \frac{\sigma}{\sqrt{2}}\),
    \begin{align*}
        & \sum_{k < j \leq \sigma} \Floor*{f(j)} = \int_{k}^{\sigma} f(\nu) \dl{\nu} - \thalf(\sigma - k) - \thalf f(k) + r \, ,\\
        \shortintertext{with}
        & r \begin{aligned}[t]
            & = \begin{aligned}[t]
                & \thalf(\sigma - k) + \thalf f(k) - f(\sigma)\sigma + k f(k) - \rho_1(f(k))k + \rho_1(f(\sigma))\sigma \\
                & {} - \thalf(\Floor{f(k)} - \Floor{f(\sigma)}) + \Floor{f(\sigma)}\Floor{\sigma} - k \Floor{f(k)} + \bigO(\sigma^{\third[2]})
            \end{aligned} \\
            & = f(k) - \Floor{f(k)} + \rho_1(\sigma)f(\sigma) + \bigO(\sigma^{\third[2]}) \\
            & = \bigO(\sigma^{\third[2]}) \, .
        \end{aligned}
    \end{align*}
    This extends to \(k \geq 0\) via splitting the sum at \(\frac{\sigma}{\sqrt{2}}\). The weighted sum then follows as in \zcref{prop:BesselStrongConvexity}.
\end{proof}

Next, a variety of integrals need to be calculated. As \(\sigma \to \infty\) then both \(\eta_1\) and \(\eta_2\) split naturally into first order terms, \(\bigO(\sigma)\), and a second order term, \(\bigO(1)\). The first order terms evaluate in terms of the incomplete Beta function. For reference recall for \(x\), \(y > 0\), and \(0 \leq z \leq 1\),
\begin{align*}
    & \BetaFDef \, , \\
    & \IncompleteBetaDef \, .
\end{align*}
\LemmaIncompleteBetaAsymptotics
\begin{proof}
    Let \(x = \sigma \sqrt{t}\) then
    \begin{align*}
        \int_{\beta}^{\sigma} x^n \sqrt{\sigma^2 - x^2} \dl{x} 
            &= \thalf \sigma^{n + 2} \int_{\frac{\beta^2}{\sigma^2}}^{1} t^{\half[n - 1]}(1 - t)^{\half} \dl{t} \\
            &= \thalf \sigma^{n + 2} \left(B\left(\thalf[n + 1], \thalf[3]\right) - B_{\frac{\beta^2}{\sigma^2}}\left(\thalf[n + 1], \thalf[3]\right)\right) \, .
    \end{align*}
    Consider the asymptotics of \(B_z(x, y)\) for fixed \(x, y\) as \(z \to 0\). Note that \((1 - t)^{y - 1} = 1 + \bigO(t)\) uniformly over compact subsets of \([0, 1)\) and so
    \begin{equation*}
        B_z(x, y) = \int_{0}^{z} t^{x - 1}(1 + \bigO(t)) \dl{t} =  \frac{z^x}{x}(1 + \bigO(z)) \, .
    \end{equation*}
    Hence, the claimed asymptotics follows with \(z = \frac{\beta^2}{\sigma^2}\), \(x = \half[n + 1]\).
\end{proof}

The second order term can be evaluated in terms of the following function. For \(p > 0\), \(q > 0\) and \(z \geq 0\) let
\begin{equation}\label{eqn:edge-integral}
    \IDef \, .
\end{equation}
\PropEdgeIntegral
\noindent The proof of \zcref{prop:edge-integral} requires the following Proposition.
\BetaDimensionalReduction
\begin{proof}[\zcref{prop:edge-integral} Proof]
    Using \zcref{eqn:BetaDimensionalReduction} and 
    \begin{equation*}
        \int_{0}^{\half[\pi]} (1 + z^2 \sin^2(\theta))^{-\half} \sin(\theta) \dl{\theta} = \frac{\arctan(z)}{z}
    \end{equation*}
    then 
    \begin{align*}
        I(z; p, 1) 
            & = 2 \int_{[0, \half[\pi]]^2} \frac{\cos^{2p - 1}(\theta) \sin^2(\theta) \sin(\phi)}{\spar*{1 + z^2 \sin^2(\theta) \sin^2(\phi)}^{\half}} \dl{\mu} \\
            & = \half \int_{[0, 1]^2} \frac{(1 - x)^{-\half} y^{\half} (1 - y)^{p - 1}}{\spar*{1 + z^2 x y}^{\half}}\dl{\mu} \\
            & = \half B(\thalf, p) \int_{0}^{1} (1 + z^2 x)^{-\half}(1 - x)^{p - \half} \dl{x} \, ,
    \end{align*}
    as required. Furthermore,
    \begin{align*}
        z^2 I(z;& p + 1, 1) \\
        & = 2 \int_{0}^{\half[\pi]} \cos(\theta)^{2 p + 1} \sin(\theta) z \arctan(z \sin(\theta)) \dl{\theta} \\
            & = \begin{aligned}[t]
                & \eval{\cos(\theta)^{2 p} \sin(\theta)^2 z \arctan(z \sin(\theta))}_0^{\half[\pi]} + p z^2 I(z; p, 2) \\
                & {} - \int_{0}^{\half[\pi]} \cos(\theta)^{2p + 1} \frac{z^2 \sin(\theta)^2}{1 + z^2 \sin(\theta)^2} \dl{\theta}
            \end{aligned} \\
            & = \begin{aligned}[t]
                & p z^2 \spar*{I(z; p, 1) - I(z; p + 1, 1)} - B(p + 1, \thalf) + \int_{0}^{\half[\pi]} \cos(\theta)^{2 p} \frac{z \cos(\theta)}{z (1 + z^2 \sin(\theta)^2)} \dl{\theta}
            \end{aligned} \\
            & = p \left((1 + z^2)I(z; p, 1) - z^2 I(z; p + 1, 1) - B(p, \thalf[3])\right) \, .
    \end{align*}
\end{proof}

The remaining terms are either estimated in an elementary fashion, or are of the form \(\int \rho_1 \dots\), which can be estimated with the following Lemma.
\LemmaSawtoothIntegralEstimate
\begin{proof}
    Denote the second Euler-Maclaurin error function as \(\rho_2\), defined as
    \begin{equation*}
        \rho_2(x) \coloneq {\half}\spar*{\rho_1(x)^2 - \frac{1}{12}}.
    \end{equation*}
    Then,
    \begin{equation*}
        \Mod*{\int_{a}^{b} \rho_1(x) f(x) \dl{x}}
        = \Mod*{\int_{a}^{b} f \dl{\rho_2}}
        = \Mod*{\eval{(\rho_2 f)}_a^b - \int_{a}^{b} \rho_2 \dl{f}}
        \leq \frac{1}{6}\spar*{ \Mod*{\eval{f}_a^b} + \int_{a}^{b} \dl{\lvert f \rvert}}
        \leq \frac{\mathcal{V}_a^b(f)}{3} \, .
    \end{equation*}
\end{proof}

These results allow proving an analogous theorem to \zcref{thm:main}, with coefficients in terms of \(n\), \(L\) and \(R\).
\TheoremTwoTermAsymptotics
\begin{proof}
    Using \zcref{prop:bessel-sum-asymptotics}, \zcref{prop:mod-bessel-sum-asymptotics}, \zcref{prop:bessel-lattice-point-func-approx}, \zcref{prop:mod-bessel-lattice-point-func-approx} and \zcref{prop:TruncateBesselSum},
    \begin{equation*}
        \N(\sigma) = \begin{aligned}[t]
            & \sum_{k = 0}^{\Floor*{\sqrt{2} \sigma \R}} \omega_{n - 2, k} + 2 \sum_{k = 0}^{\Floor*{\sigma \R - {\half} (\sigma \R)^{\gamma}}} \omega_{n - 2, k} \Floor*{\frac{\eta_1(k + \beta, \sigma \R) - \tquarter[\pi]}{\pi} + \fcall{\bigO}{\sigma^{\half[1 - 3 \gamma]}}}\\
            & {} + \sum_{k = 0}^{\Floor*{\sigma \R}} \omega_{n - 2, k} \Floor*{\frac{\eta_2(k + \beta, \sigma \R + \beta)}{\pi} + \fcall{\bigO}{\sigma^{-1}}} \\
            & {} + \sum_{k = 0}^{\Floor*{\sigma \R}} \omega_{n - 2, k} \Floor*{\frac{\eta_2(k + \beta, \sigma \R + \beta) + \thalf[\pi]}{\pi} + \fcall{\bigO}{\sigma^{-1}}} \\
            & {} + \bigO(\sigma^{\half[5 \gamma] - \half} + \sigma^\gamma) \sigma^{\n - 3} \, .
        \end{aligned}
    \end{equation*}
    Considering each term, first
    \begin{equation*}
        \sum_{k = 0}^{\Floor*{\sqrt{2} \sigma \R}} \omega_{n - 2, k} = \frac{2^{\half[\n]}}{(\n - 2)!}(\sigma \R)^{\n - 2} + \bigO(1) \, .
    \end{equation*}
    Next, note that
    \begin{equation*}
        \binom{\Floor*{x} + \n - 4}{\n - 4} = \frac{(\Floor*{x} + 1)_{\n - 3}}{(\n - 3)!} = \frac{(x + {\thalf})_{\n - 3} + \rho_1(x)r(x)}{(\n - 3)!} \, ,
    \end{equation*}
    for some polynomial \(r\) of max degree \(n - 4\). Hence,
    \begin{align*}
        & \sum_{k = 0}^{\Floor*{\sigma \R - {\half} (\sigma \R)^{\gamma}}} \omega_{n - 2, k} \Floor*{\frac{\eta_1(k + \beta, \sigma \R) - \tquarter[\pi]}{\pi} + \fcall{\bigO}{\sigma^{\half[1 - 3 \gamma]}}} \\
        & \quad = \sum_{k = 0}^{\Floor*{\sigma \R - {\half} (\sigma \R)^{\gamma}}} \binom{k + \n - 4}{\n - 4} \sum_{j = 0}^{\Floor*{\sigma \R - {\half} (\sigma \R)^{\gamma}} - k} \omega_{1, j} \Floor*{\frac{\eta_1(j + \beta, \sigma \R) - \tquarter[\pi]}{\pi} + \fcall{\bigO}{\sigma^{\half[1 - 3 \gamma]}}} 
        \shortintertext{\rightline{\zcref{thm:HarmonicWeightSlicing}}}
        & \quad = \begin{aligned}[t]
            & \sum_{k = 0}^{\Floor*{\sigma \R - {\half} (\sigma \R)^{\gamma}}} \binom{k + \n - 4}{\n - 4} \spar*{2 \int_{k}^{\sigma \R - {\half} (\sigma \R)^{\gamma}} \frac{\eta_1(\nu + \beta, \sigma \R)}{\pi} - {\tquarter} \dl{\nu} - (\sigma \R - k)} \\
            & {} + \bigO(\sigma^{\third[2]} + \sigma^{\half[3](1 - \gamma)} + \sigma^\gamma) \sigma^{\n - 3}
        \end{aligned}
        \shortintertext{\rightline{\zcref{prop:BesselStrongConvexity}}}
        & \quad = \begin{aligned}[t]
            & 2 \int_{\beta}^{\sigma \R} \frac{(\Floor*{\nu - \beta} + 1)_{\n - 3}}{(\n - 3)!} \frac{\eta_1(\nu, \sigma \R) - \tquarter[\pi]}{\pi} \dl{\nu} - \frac{(\sigma \R)^{\n - 2}}{(\n - 2)!} + \bigO(\sigma^{\third[2]} + \sigma^{\half[3](1 - \gamma)} + \sigma^\gamma) \sigma^{\n - 3}
        \end{aligned} \\
        & \quad = \begin{aligned}[t]
            & \frac{B(\thalf[\n - 2], \thalf[3])}{\pi (\n - 1)(\n - 3)!} (\sigma \R)^{\n - 1} - \spar*{\frac{1}{2 (\n - 2)!} + \frac{I(1; \thalf[\n - 2], 1)}{\pi (\n - 3)!}} (\sigma \R)^{\n - 2} \\
            & {} + \bigO(\sigma^{\third[2]} + \sigma^{\half[3](1 - \gamma)} + \sigma^\gamma) \sigma^{\n - 3} \, .
        \end{aligned}
        \shortintertext{\rightline{\zcref{lemma:SawtoothIntegralEstimate,lemma:IncompleteBetaAsymptotics}}}
    \end{align*}
    Similarly for \(\eta_2\) but using \zcref{prop:ModifiedBesselStrongConcavity}
    \begin{align*}
        & \sum_{k = 0}^{\Floor*{\sigma \R}} \omega_{n - 2, k} \Floor*{\frac{\eta_2(k + \beta \R, (\sigma + \beta) \R)}{\pi} + \fcall{\bigO}{\sigma^{-1}}} \\
        & \begin{aligned}[t]
            & = \sum_{k = 0}^{\Floor{\sigma \R}} \binom{k + \n - 4}{\n - 4} \sum_{j = 0}^{\Floor*{\sigma \R} - k} \omega_{1, j} \Floor*{\frac{\eta_2(k + \beta, \sigma \R + \beta)}{\pi} + \fcall{\bigO}{\sigma^{-1}}} \\
            & = \sum_{k = 0}^{\Floor{\sigma \R}} \binom{k + \n - 4}{\n - 4} \spar*{2 \int_{k}^{\sigma \R} \frac{\eta_2(\nu + \beta, \sigma \R + \beta)}{\pi} \dl{\nu} - (\sigma \R - k) + \bigO(\sigma^{\third[2]})} \\
            & = 2 \int_{0}^{\sigma \R} \frac{(\Floor*{\nu} + 1)_{\n - 3}}{(\n - 3)!} \frac{\eta_2(\nu + \beta, \sigma \R + \beta)}{\pi} \dl{\nu} - \frac{(\sigma \R)^{\n - 2}}{(\n - 2)!} + \bigO(\sigma^{\third[2]})\sigma^{\n - 3} \\
            & = \begin{aligned}[t]
                & \frac{B(\thalf[\n - 2], \thalf[3])}{\pi (\n - 3)!}  \frac{\oL}{\R} (\sigma \R)^{\n - 1} + 2 \Biggl(\beta^2 B(\thalf[\n - 3], \thalf[3]) + \thalf[\beta]\spar*{B(\thalf[\n - 2], {\thalf}) - B(\thalf[\n - 1], {\thalf})} \\
                & {} + {\tquarter}B(\thalf[\n - 1], \thalf[3]) \Biggr)  \frac{\oL}{\R} \frac{(\sigma \R)^{\n - 2}}{\pi (\n - 3)!} - \frac{I(1; \thalf[\n - 2], 1)}{\pi (\n - 3)!}(\sigma \R)^{\n - 2} - \frac{(\sigma \R)^{\n - 2}}{(\n - 2)!} + \bigO(\sigma^{\third[2]}) \sigma^{\n - 3} \, .
            \end{aligned}
        \end{aligned}
    \end{align*}
    Note that
    \begin{equation*}
        \beta^2 B(\thalf[\n - 3], \thalf[3]) + \thalf[\beta]\spar*{B(\thalf[\n - 2], {\thalf}) - B(\thalf[\n - 1], {\thalf})}  + {\tquarter}B(\thalf[\n - 1], \thalf[3]) = \frac{\n - 2}{4} B(\thalf[\n], {\thalf}) \, ,
    \end{equation*}
    then it is immediate that for the last sum
    \begin{multline*}
        \sum_{k = 0}^{\Floor*{\sigma \R}} \omega_{n - 2, k} \Floor*{\frac{\eta_2(k + \beta, \sigma \R + \beta) + \half[\pi]}{\pi} + \fcall{\bigO}{\sigma^{-1}}} = \\ \frac{B(\thalf[\n - 2], \thalf[3])}{\pi (\n - 3)!}  \frac{\oL}{\R}(\sigma \R)^{\n - 1} + \frac{(\n - 2)B(\thalf[\n], {\thalf})}{2 \pi (\n - 3)!}  \frac{\oL}{\R} (\sigma \R)^{\n - 2} - \frac{I(1; \thalf[\n - 2], 1)}{\pi (\n - 3)!}(\sigma \R)^{\n - 2} + \bigO(\sigma^{\third[2]}) \sigma^{\n - 3} \, .
    \end{multline*}
    Collecting terms
    \begin{align*}
        \N(\sigma) = \begin{aligned}[t]
            & \spar*{\frac{2 B(\thalf[\n], \thalf)}{\pi (\n - 1)!} + \frac{2 B(\thalf[\n], {\thalf})}{\pi (\n - 2)!}   \frac{\oL}{\R}} (\sigma \R)^{\n - 1} \\
            & {} + \spar*{\frac{(\n - 2)B(\thalf[\n], \thalf)}{\pi (\n - 3)!}  \frac{\oL}{\R} + \spar*{2\frac{2^{\half[\n]} - 1}{(\n - 2)!} + \frac{4 I(1; \thalf[\n - 2], 1)}{\pi (\n - 3)!} }} (\sigma \R)^{\n - 2}  \\
            & {} + \bigO(\sigma^{\third[2]} + \sigma^{\half[3](1 - \gamma)} + \sigma^\gamma + \sigma^{\half[5 \gamma - 1]}) \sigma^{\n - 3} \, .
        \end{aligned}
    \end{align*}
    The optimal choice is \(\gamma = \half\) leading to an error of \(\sigma^{n - 2 - \quarter}\) as claimed.
\end{proof}

\zcref{thm:two-term-asymptotics} can now be used to show \zcref{thm:main} with the following Lemma.
\LemmaBetaBallVolumeRelations
\begin{proof}
    It is well known that
    \begin{equation*}
        \measure*{\mathbb{B}_p} = \prod_{i = 1}^{p} B(\thalf[i] + \thalf, \thalf) = \frac{\Gamma(\thalf)^p}{\Gamma(\thalf[p] + 1)} \, .
    \end{equation*}
    Hence,
    \begin{align*}
        \measure*{\mathbb{B}_p}^2 
            & = \frac{\pi^p}{\Gamma(\thalf[p] + 1)^2} \\
            & = (2 \pi)^p \frac{\Gamma(\thalf)}{2^p \Gamma(\thalf[p] + \thalf) \Gamma(\thalf[p] + 1)} \frac{\Gamma(\thalf[p] + \thalf)\Gamma(\thalf)}{\pi \Gamma(\thalf[p] + 1)} \\
            & = \frac{(2 \pi)^p}{p!} \frac{B(\thalf[p] + \thalf, \thalf)}{\pi}\, , \\
    \end{align*}
    and
    \begin{equation*}
        \measure*{\mathbb{B}_p} \measure*{\mathbb{B}_{p - 1}} = \frac{\Gamma(\thalf)^{2p - 2} \Gamma(\thalf)}{\Gamma(\thalf[p] + \thalf)\Gamma(\thalf[p] + 1)} = \frac{2 (2 \pi)^{p - 1}}{p!} \, .
    \end{equation*}
\end{proof}

\begin{proof}[Proof of \zcref{thm:main}]
Recall that
\begin{align*}
    & \measure*{\partial\Omega} = 2 \measure*{\mathbb{B}_{n - 1}} R^{n - 1} + 2 L (n - 1) \measure*{\mathbb{B}_{n - 1}} R^{n - 2} \, , \\
    & \measure*{\partial^2 \Omega} = 2 (n - 1) \measure*{\mathbb{B}_{n - 1}} R^{n - 2} \, ,\\
    & \int_{\partial\Omega} \kappa = \frac{n - 2}{(n - 1) R} 2 L (n - 1) \measure*{\mathbb{B}_{n - 1}} R^{n - 2} \, .
\end{align*}
Hence, using \zcref{lemma:BetaBallVolumeRelations}
\begin{align*}
    & \begin{aligned}[t]
        \frac{\measure*{\mathbb{B}_{n -1}}}{(2 \pi)^{n - 1}} \measure*{\partial\Omega} & = \frac{\measure*{\mathbb{B}_{n -1}}^2}{(2 \pi)^{n - 1}} \left(2 R^{n - 1} + 2 L (n - 1) R^{n - 2}\right) \\
        & = \frac{B(\thalf[n], \thalf)}{\pi(n - 1)!}\left(2 R^{n - 1} + 2 L (n - 1) R^{n - 2}\right) \, ,
    \end{aligned} \\
    & \begin{aligned}[t]
        & \frac{\measure*{\mathbb{B}_{n -1}}}{(2 \pi)^{n - 1}} \frac{(n - 2)(n - 1)}{2} \int_{\partial \Omega} \kappa \\
        & \qquad = \frac{\measure*{\mathbb{B}_{n - 1}}^2}{(2 \pi)^{n - 1}} \left((n - 2)^2(n - 1)\right) L R^{n - 3} \\
        & \qquad = \frac{n - 2}{\pi (n - 3)!} B(\thalf[n], \thalf) L R^{n - 3} \, ,
    \end{aligned} \\
    & \begin{aligned}[t]
        \frac{\measure*{\mathbb{B}_{n - 2}}}{(2 \pi)^{n - 2}} \measure*{\partial^2 \Omega} & = 2 (n - 1) \frac{\measure*{\mathbb{B}_{n - 2}}\measure*{\mathbb{B}_{n - 1}}}{(2 \pi)^{n - 2}}R^{n - 2} = \frac{4 R^{n - 2}}{(n - 2)!} \, ,
    \end{aligned} \\
    & \begin{aligned}[t]
        \frac{\measure*{\mathbb{B}_{n - 2}}}{(2 \pi)^{n - 2}}\measure*{\partial^2 \Omega} G'_{n - 1, 1} 
        & = \frac{4 R^{n - 2}}{\pi (n - 3)!} I(1; \thalf[n - 2], 1) \,.
    \end{aligned}
\end{align*}
\end{proof}

\appendix

\section{Equivalence of Edge Term In \texorpdfstring{
    \cite[Theorem 1.1]{SteklovCuboids}}{
        [The Steklov Spectrum of Cuboids, Theorem 1.1]}}
\label{section:edge-equivalence}
In \cite{SteklovCuboids} it is shown that the spectral counting function of the Steklov spectrum of a cuboid, \(\Omega\), of dimension \(n \geq 3\), admits the following asymptotic expansion as \(\sigma \to \infty\)
\begin{equation*}
    N(\sigma) = \frac{\Mod*{\mathbb{B}_{n - 1}}}{(2 \pi)^{n-1}} \Mod*{\partial \Omega} + \frac{\Mod*{\mathbb{B}_{n - 2}}}{(2 \pi)^{n-2}}\spar*{2^{\thalf[n] - 1} - {\half}} - \frac{2}{\pi^{n - 1}}G_{n - 1, 1} + \bigO\spar*{\sigma^{d - 2 - \min\spar*{{\third}, (d - 1)^{-1}}}} \,,
\end{equation*}
where
\begin{equation*}
    G_{n, 1} \coloneq \int_{[0, \half[\pi]]^{n - 1}} \fcall{\arccot}{\prod_{i = 1}^{n - 1}\csc(\theta_i)} \prod_{i = 1}^{n - 1} \sin^i(\theta_i) \dl{\mu}(\theta_1, \dots, \theta_{n - 1}) \,.
\end{equation*}
Comparing with the asymptotic expansion from \zcref{thm:main} then the two expansions would be consistent, as cuboids have no curvature, if
\begin{equation*}
    \frac{2 G_{n, 1}}{\pi^n} = \frac{\lvert \mathbb{B}_{n-1} \rvert}{(2 \pi)^{n- 1}} G'_{n, 1}
\end{equation*}
where
\begin{equation*}
    G'_{n, 1} \coloneq \frac{\int_{0}^{1} (1 + x)^{-\half}(1 - x)^{\half[n] - 1} \dl{x}}{B(\thalf, \thalf[n])} \,.
\end{equation*}
Indeed, this can be shown using \zcref{eqn:BetaDimensionalReduction} as follows
\begin{align*}
    G_{n, 1} 
        & = \int_{[0, \half[\pi]]^{n - 1}} \fcall{\arccot}{\prod_{i = 1}^{n - 1}\csc(\theta_i)} \prod_{i = 1}^{n - 1} \sin^i(\theta_i) \dl{\mu} \\
        & = \int_{[0, \half[\pi]]^{n - 1}} \fcall{\arctan}{\prod_{i = 1}^{n - 1}\sin(\theta_i)} \prod_{i = 1}^{n - 1} \sin^i(\theta_i) \dl{\mu} \\
        & = \int_{[0, \half[\pi]]^{n}} \frac{\prod_{i = 1}^{n} \sin^i(\theta_i)}{\left(1 + \prod_{i = 1}^n \sin^2(\theta_i)\right)^{\half}} \dl{\mu} \\
        & = \frac{1}{2^n} \int_{[0, 1]^{n}} \left(1 + \prod_{i = 1}^{n} x_i\right)^{-\half}\prod_{i = 1}^n \frac{x_i^{\half[i - 1]}}{(1 - x^i)^{\half}} \dl{\mu} \\
        & = \frac{\Gamma(\thalf)^n}{2^n \Gamma(\thalf[n])} \int_{0}^{1} (1 + x)^{-\half}(1 - x)^{\half[n] - 1} \dl{x} \,.
\end{align*}

Furthermore, in \cite{SteklovCuboids} only \(G_{2, 1}\) and \(G_{3, 1}\) were calculated explicitly. Noting that 
\begin{equation*}
    I(1; \thalf[n - 1], 1) = \frac{\pi}{n - 1}G'_{n, 1}
\end{equation*}
then the recurrence relation from \zcref{prop:edge-integral} passes through to the following recurrence relation for \(G'_{n, 1}\)
\begin{equation}
    G'_{0, 1} = \half[\sqrt{2}], G'_{1, 1} = \half, G'_{n + 2, 1} = 2\spar*{G'_{n, 1} - \frac{1}{n B(\thalf[n], \thalf)}} \,,
\end{equation}
allowing efficient calculation of \(G'_{n, 1}\) for all \(n\), and hence, also for \(G_{n, 1}\).

\section{Proof of \texorpdfstring{\zcref{prop:BetaDimensionalReduction}}{Integral Dimension Reduction Proposition}}
For completeness, a proof of \zcref{prop:BetaDimensionalReduction} is given.
\begin{proof}
    First consider \((x, y) \in [0, 1]^2\) with the standard Lebesgue measure. The following is a bijection outside a set of measure zero, \(x' = xy, y'(1 - x') = y(1 - x)\). The corresponding change in measure is
    \begin{align*}
        y \dl{\mu} = (1 - x') \dl{\mu'} \,.
    \end{align*}
    Hence, given \(a, b > 0\) and \(f\) such that \(f(x)(1 - x)^{a + b - 1} \in L^1([0, 1])\) then
    \begin{align*}
        & \int_{[0, 1]^2} f(x y) (1 - x)^{a - 1} y^{a} (1 - y)^{b - 1} \dl{\mu} \\
        & \qquad \begin{aligned}[t]
            & = \int_{[0, 1]^2} f(x') (1 - x')^{a + b - 1} y'^{a - 1} (1 - y')^{b - 1} \dl{\mu'} \\
            & = B(a, b) \int_{0}^{1} f(x') (1 - x')^{a + b - 1} \dl{x'} \, .
        \end{aligned}
    \end{align*}
    The result then follows by induction.
\end{proof}

\printbibliography

{
    \renewcommand\thefootnote{}
    \footnotetext[0]{\parbox[t]{\textwidth - 15pt}{
        Address: Spencer Bullent, Department of Mathematics, University College London, London, WC1H 0AY, UK
        Email: \href{emailto:spencer.bullent.23@ucl.ac.uk}{spencer.bullent.23@ucl.ac.uk}}
    }
}

\end{document}